# Reduced Order Modeling based Inexact FETI-DP solver for lattice structures

T. Hirschler[a,b,*], R. Bouclier[c,d], P. Antolin[b], A. Buffa[b,e]

[a]*ICB, UMR6303 CNRS, Université de Technologie de Belfort-Montbéliard, Sévenans, France*
[b]*Institute of Mathematics - Chair of Numerical Modelling and Simulation*
*École Polytechnique Fédérale de Lausanne, Switzerland*
[c]*IMT, Université de Toulouse, UPS-UT1-UT2-INSA-CNRS, Toulouse, France*
[d]*ICA, Université de Toulouse, INSA-ISAE-Mines Albi-UPS-CNRS, Toulouse, France*
[e]*Istituto di Matematica Applicata e Tecnologie Informatiche "Enrico Magenes"*
*Consiglio Nazionale delle Ricerche, Pavia, Italy*

**Abstract**

This paper addresses the overwhelming computational resources needed with standard numerical approaches to simulate architected materials. Those multiscale heterogeneous lattice structures gain intensive interest in conjunction with the improvement of additive manufacturing as they offer, among many others, excellent stiffness-to-weight ratios. We develop here a dedicated HPC solver that benefits from the specific nature of the underlying problem in order to drastically reduce the computational costs (memory and time) for the full fine-scale analysis of lattice structures. Our purpose is to take advantage of the natural domain decomposition into cells and, even more importantly, of the geometrical and mechanical similarities among cells. Our solver consists in a *so-called* inexact FETI-DP method where the local, cell-wise operators and solutions are approximated with reduced order modeling techniques. Instead of considering independently every cell, we end up with only few principal local problems to solve and make use of the corresponding principal cell-wise operators to approximate all the others. It results in a scalable algorithm that saves numerous local factorizations. Our solver is applied for the isogeometric analysis of lattices built by spline composition, which offers the opportunity to compute the reduced basis with macro-scale data, thereby making our method also multiscale and matrix-free. The solver is tested against various 2D and 3D analyses. It shows major gains with respect to black-box solvers; in particular, problems of several millions of degrees of freedom can be solved with a simple computer within few minutes.

## 1. Introduction

Additive Manufacturing (AM) and especially its metal variants constitute today a reality for the fabrication of industrially-relevant high-performance parts and products [6]. In this context, AM has pushed the development of architected materials to a new stage. These can be viewed as specific multiscale heterogeneous lattice structures that are inspired, in particular, from biology such as bone architectures. In practice, they usually consist in repeating a pattern, the unit-cell, made of a combination of struts connected in different ways, all over a macro shape (see Figure 1). The cell topology, geometry and material can be tailored for achieving specific performance requirements [51]. Starting with unprecedented weight saving while maintaining the stiffness and strength of the structure, lattices can also accommodate novel behaviors (*e.g.*, highly stretchable and auxetic) and multi-functionality, which make them attractive for lightweight components, energy absorbing, thermal management, design of medical implants, to name a few [5, 58].

With the improvement of AM, representative lattices may now come with large numbers of cells of complex topology and shape, which make their numerical simulation very challenging. In particular, the computational cost and memory requirements become quickly tremendous or even intractable if standard methods are blindly applied. Thus, the efficient prediction of the structural response of lattice structures is becoming a hot topic nowadays.

---

[*]Corresponding author
  *Email address:* `thibaut.hirschler@utbm.fr` (T. Hirschler)

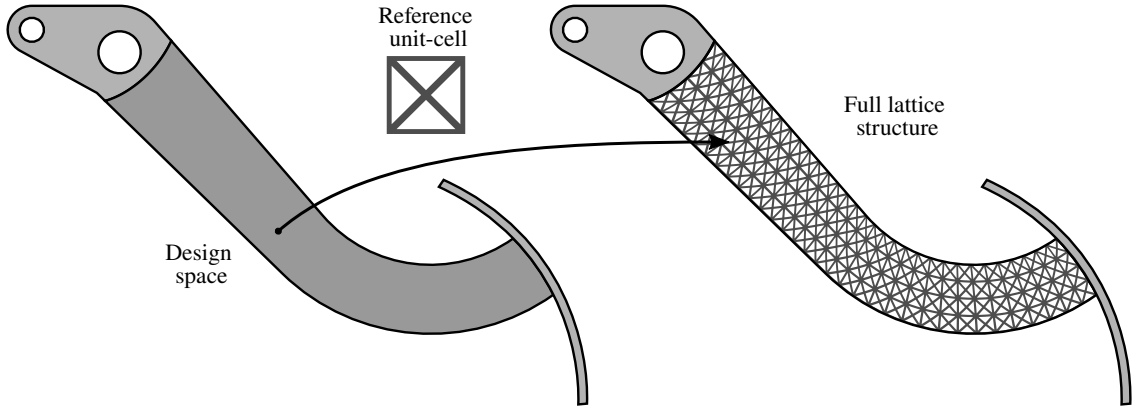

Figure 1: Example of a lattice structure (inspired from [49]): The pedal crank is composed of a lattice structure (right) which allows drastic weight saving compared to a classic design (left).

As a solution, one may be tended to draw inspiration from the multitude multiscale methods based on homogenization that have been developed over the years for multiscale heterogeneous structures. The goal of these approaches is to compute a homogenized macro-scale behavior from its underlying heterogeneous microstructure by different interchange of information between scales. Among others, let us mention the multilevel FEM (FE$^2$ [4, 19, 53, 74]) which computes homogenized fields from the micro-scale and assigns them at each integration point of the macro-scale model; the Multiscale FEM (MsFEM [30, 75]) where numerically computed basis functions encode the micro-scale heterogeneity; global/local coupling [70] that allows to replace the homogenized global model by the fine high-fidelity model in a certain local region; or direct numerical homogenization [22, 34, 47, 64] where macro-scale material parameters are identified. However, these methods mostly rely on a true separation of scales to alleviate size effects [50, 76]. This is usually not the case with lattice structures since 3D printers notwithstanding constraint the achievable length scale range of the cells. Furthermore, the macro-shape often follows a rather slender geometry (plate, shell), thereby exhibiting few cells in one direction. Alternatively, specific structural models such as sophisticated 3D beams [71, 72] could be used to model the full architecture of the lattice while reducing the complexity of the problem. Yet, such methods may not be suitable in case of small cells (thus involving thick struts) and cannot be applied to more general lattices, (*e.g.*, skeletal-TPSM or sheet-TPMS based lattices [5]). As a result, it seems that instead of acting at the modeling stage, it would be desirable to work on the solution stage for lattices; that is, to derive HPC solvers capable of solving the high-fidelity, fine-scale problem efficiently.

In this spirit, immersed domain methods, such as the Finite Cell Method [52, 59], have recently been applied for the full fine-scale analysis of lattice structures [40–42]. Such approaches switch the effort of mesh generation to the proper integration of cut elements and the sophisticated treatment of conditioning problems. Therefore, when used in combination with image-based modeling [14, 57, 66], they make it easier to take into account the (strut-level) geometrical defects due to the AM process (see [40–42] again). HPC extensions of such technologies exist [2, 33]. They are generic in terms of applications but they may not appear optimal in the specific context of lattice structures. In particular, they do not make use of the geometric similarity among different cells in the numerical solution, what allows a drastic reduction in terms of memory and computational cost, as demonstrated for the assembly of the numerical operators in [28]. The objective of this work is to bring forward this idea of benefiting from the repetitive character of the lattice geometries to build a dedicated HPC solver.

In order to do so, our approach relies on three main ingredients. Firstly, it is based on the Dual-Primal variant of the Finite Element Tearing and Interconnecting method (FETI-DP), that, after its appearance in [18] has been widely developed and adopted in high-performance computational structural mechanics (see [35, 38, 63, 65] to name a few). This algorithm belongs to the class of non-overlapping Domain Decomposition (DD) solvers that appear highly relevant here given the natural splitting of the lattice structures into several subdomains (each cell can be associated to a subdomain). DD methods are preconditioned iterative algorithms where, at each iteration, subdomain-independent local systems are solved in parallel. The interface conditions are recovered through the iterative process. To obtain a scalable algorithm, a small coarse problem also needs to be solved at each iteration. We note that DD solvers have been used in the past to successfully compute multiscale heterogeneous structures (see, *e.g.*, [15, 43, 44]). In FETI-DP, the continuity constraints on the displacement at the subdomain corners are forced to hold throughout the iterative procedure, which naturally leads to solving a coarse problem at each iteration, while the other constraints are enforced by the use of Lagrange multipliers. The second ingredient of our approach is to consider the framework of inexact FETI-DP algorithms [36, 55]. In particular, we develop an inexact version that offers the opportunity to avoid solving numerous local (subdomain-wise) systems. For this, we iterate on the initial complete saddle point problem and design a block preconditioner consistent with the mathematical background established in [48, 60]. Thirdly, we investigate the repetitiveness of the cells by applying a Reduced Order Modeling (ROM) approach [25, 54] (especially, a greedy technique) that extracts the "principal" cells in terms of stiffness. The resulting few local stiffness operators are then used in the preconditioner of the inexact FETI-DP to build reduced bases for the efficient approximate solutions of the numerous local systems. It results in a scalable algorithm that saves numerous local factorizations and thus offers the opportunity to deeply reduce the computational time with respect to the standard FETI-DP solver, in particular when the macro mapping is simple which is usually the case in real lattice structures. In this contribution, we place ourselves at the design stage of the product; that is, we only consider perfect lattices. Yet, we believe that the proposed framework owns all the necessary ingredients to be easily extended to treat lattices with geometrical imperfections due to AM. In addition, we assume here, as a starting point, that the lattice is linear elastic (at the scale of the cell-struts) and we place ourselves in the context of small displacements and rotations.

Furthermore, although not mandatory in our methodology, we adopt the Computer-Aided Design (CAD) paradigm based on spline composition [16] to generate the lattice geometries. Such a geometrical modeling offers great flexibility to design lattices since the local (cell) and global (macro) geometries are naturally parameterized by splines. In addition, it is fully consistent with IsoGeometric Analysis (IGA) [12, 31], thus allowing to directly analyze such models (see, *e.g.*, [1, 27]) while avoiding possibly delicate meshing procedures. A first step toward the efficient analysis of lattices modeled by spline composition has been recently performed with the development of a fast multiscale assembly procedure for the IGA operators (see previous contribution [28]). The strategy makes use of look-up tables with precomputed integrals associated to the cell pattern and macro-fields that encode the mechanical behavior related to the macro-shape. Benefiting from this data structure, we are actually able to extract the "principal" cells in our solver by acting only on the macro-fields, thereby avoiding forming all the local stiffness operators. In this context, our ROM-based inexact FETI-DP solver thus incorporates additional attractive features: It is multiscale and tends to be matrix-free. From this point of view, the present work also further advances the state-of-the-art of DD-IGA that represents an active research area (see [3, 39] for the origins and then [7, 26, 29, 73] to name a few).

The remainder of this paper is organized as follows. In Section 2, we briefly review the FETI-DP method starting from the complete saddle point formulation and give the key ingredients to build efficient block preconditioners which is the basis of our solver. Then, we describe in Section 3 the developed ROM-based inexact strategy to benefit from the repetitive character of lattice geometries in the FETI-DP solver. Making use of algebraic equations, the presentation is kept general at this stage to highlight the applicability of our solver to any discretization scheme. In Section 4, we specify the application of the methodology in case of a CAD modeling based on spline composition, which results in an isogeometric solver that is also multiscale and quasi matrix-free. Next, numerical results to assess the efficiency of our solver both in terms of memory and computational cost reduction are reported in Section 5 and, finally, Section 6 contains some conclusions and future directions of research.

## 2. FETI-DP and saddle point problems

### 2.1. Domain decomposition

A lattice structure is by construction made of multiple cells that are packed together and fill a particular area of a macro-geometry (recall Figure 1). In other words, the domain $\Omega$ describing a lattice structure can be naturally seen as a non-overlapping collection of subdomains $\Omega^{(s)}$, each of them being associated to a lattice unit-cell. Thus, we define:

$$\bar{\Omega} = \bigcup_{s=1}^{N} \bar{\Omega}^{(s)}, \quad \text{with} \quad \Omega^{(s_1)} \cap \Omega^{(s_2)} = \emptyset \text{ if } s_1 \neq s_2, \tag{1}$$

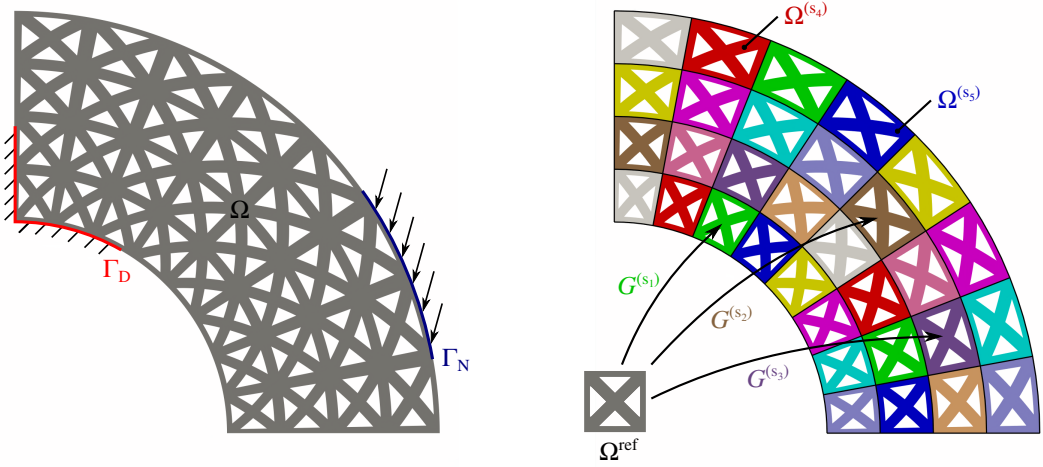

Figure 2: Problem to be solved along with its domain decomposition. The finite element spaces associated to the different subdomains are constructed from the basis functions of the reference unit-cell.

and where $N$ is the total number of cells. Let us also introduce the subdomain interfaces:

$$\Gamma = \bigcup_{s_1 \neq s_2} \partial\Omega^{(s_1)} \cap \partial\Omega^{(s_2)}. \tag{2}$$

As mentioned in the introduction, there exist specific solvers which are particularly meant to take advantage of such a domain decomposition: the so-called non-overlapping DD algorithms. The basis of this work is thus to resort to the DD technology to build an efficient solver for lattice structures.

**Remark 1.** *Note that in the general field of DD solvers, obtaining a domain decomposition of the form* (1) *actually requires a pre-processing step that can be rather complex and may involve sophisticated tools (see, e.g., [13, 24, 69]). However, here, the geometrical structure of the domain leads in a straightforward manner to a non-overlapping decomposition, making the application of DD solvers a natural choice.*

### 2.2. Dual-primal FETI method

More precisely, we investigate in this work a class of DD solvers known as FETI-DP [18], where the starting point is the construction of a particular working space in which we will iteratively look for our (finite element) solution. Then, a specific algorithm is derived which takes advantage of the particular structure of this working space.

#### 2.2.1. Model problem

We consider here the linear elasticity problem in $\Omega$ which is naturally decomposed in $N$ subdomains as expressed in Equation (1) (see Figure 2 for illustration). In each subdomain, the displacement field $u^{(s)} : \bar{\Omega}^{(s)} \to \mathbb{R}^d$, $d$ being the dimension of the problem, satisfies the following equations:

$$\text{div}(\sigma^{(s)}) + b_{\text{ext}} = 0 \quad \text{in} \quad \Omega^{(s)}, \tag{3a}$$

$$\sigma^{(s)} = C^{(s)} \varepsilon(u^{(s)}) \quad \text{in} \quad \Omega^{(s)}, \tag{3b}$$

$$u^{(s)} = u_0 \quad \text{over} \quad \Gamma_D \cap \partial\Omega^{(s)}, \tag{3c}$$

$$\sigma^{(s)} n^{(s)} = t_{\text{ext}} \quad \text{over} \quad \Gamma_N \cap \partial\Omega^{(s)}, \tag{3d}$$

where $\varepsilon$ denotes the linearized Green-Lagrange strain tensor, $\sigma$ the linearized Cauchy stress tensor, C the Hooke material tensor, and $n^{(s)}$ the outward unit normal to $\partial\Omega^{(s)}$. The external loads are of two kinds: $b_{\text{ext}}$ denotes body forces, and $t_{\text{ext}}$ stands for tractions prescribed at the boundary $\Gamma_N$ of the body. Imposed displacements $u_0$ over the boundary $\Gamma_D$ are also considered (such that the total problem is well posed). In order to recover the linear elasticity problem posed over the full domain $\Omega$, kinematic compatibility and equilibrium of the tractions along the subdomain interfaces $\Gamma$ should be added. For every couple $(s_1, s_2)$ of neighboring subdomains, this leads to write:

$$u^{(s_1)} - u^{(s_2)} = 0 \quad \text{over} \quad \Gamma^{(s_1,s_2)}, \tag{4a}$$

$$\sigma^{(s_1)} n^{(s_1)} + \sigma^{(s_2)} n^{(s_2)} = 0 \quad \text{over} \quad \Gamma^{(s_1,s_2)}, \tag{4b}$$

with $\Gamma^{(s_1,s_2)} = \partial\Omega^{(s_1)} \cap \partial\Omega^{(s_2)} \subset \Gamma$, $s_1 < s_2$, and $\|\Gamma^{(s_1,s_2)}\| > 0$.

This coupled problem is solved approximately using a finite element approach. To this purpose, let us recall the variational forms that come with the problem defined in Equation (3):

$$a^{(s)}(u, v) = \int_{\Omega^{(s)}} \sigma^{(s)}(u) : \varepsilon(v) \, dx, \tag{5a}$$

$$b^{(s)}(v) = \int_{\Omega^{(s)}} b_{\text{ext}} \cdot v \, dx + \int_{\Gamma_N \cap \partial\Omega^{(s)}} t_{\text{ext}} \cdot v \, dl. \tag{5b}$$

At this stage, let us make the following important assumption: We consider that the finite element spaces defined over the subdomains (denoted as $\mathcal{W}^{h(s)}$ and defined such that $\mathcal{W}^{h(s)} \subset \left[H^1(\Omega^{(s)})\right]^d$) are all built starting from the same functions. More specifically, we introduce a set of basis functions $N_i : \Omega^{\text{ref}} \to \mathbb{R}, i = 1, \ldots, n$, where $\Omega^{\text{ref}}$ denotes the domain associated to the reference unit-cell, and we assume that for each subdomain there exits a smooth mapping $G^{(s)} : \Omega^{\text{ref}} \to \Omega^{(s)}$ with a smooth inverse that maps each point $\boldsymbol{\xi}$ in $\Omega^{\text{ref}}$ to a point $\mathbf{x}$ in $\Omega^{(s)}$ (see Figure 2(right)). Thus, we expand the local finite element spaces in terms of basis functions, as:

$$\mathcal{W}^{h(s)} = \text{span}\left\{N_i \circ G^{(s)-1}\right\}_{i=1,\ldots,n}. \tag{6}$$

Finally, let us decompose these local spaces as:

$$\mathcal{W}^{h(s)} = \mathcal{W}_\Gamma^{h(s)} \oplus \mathcal{W}_I^{h(s)}, \tag{7}$$

where $\mathcal{W}_\Gamma^{h(s)}$ is the finite element trace space (*i.e.*, made of the basis functions that do not vanish at the subdomain boundary) and $\mathcal{W}_I^{h(s)}$ is the finite element space associated to the interior part (*i.e.*, $\forall v \in \mathcal{W}_I^{h(s)}, v_{|\partial\Omega^{(s)}} = 0$).

2.2.2. Dual-Primal continuity constraints

The local spaces $\mathcal{W}^{h(s)}$ can be combined to define a working global space for the total solution:

$$\mathcal{W}^h = \bigoplus_{s=1}^N \mathcal{W}^{h(s)} = \mathcal{W}_I^h \oplus \mathcal{W}_\Gamma^h, \quad \text{with} \quad \mathcal{W}_I^h = \bigoplus_{s=1}^N \mathcal{W}_I^{h(s)}, \mathcal{W}_\Gamma^h = \bigoplus_{s=1}^N \mathcal{W}_\Gamma^{h(s)}. \tag{8}$$

However, the members of $\mathcal{W}^h$ are in general not continuous across the subdomain interfaces. Thus, an alternative choice would be to look for solutions in the following subset:

$$\hat{\mathcal{W}}^h = \left\{u^h \in \mathcal{W}^h : u^h \text{ is continuous across } \Gamma\right\}. \tag{9}$$

The standard FETI method [17] iterates in space $\mathcal{W}^h$, whereas the alternative Balancing Domain Decomposition (BDD) method [45] iterates in space $\hat{\mathcal{W}}^h$. Within FETI-DP approaches [18], the main idea consists in iterating in an intermediate space between $\mathcal{W}^h$ and $\hat{\mathcal{W}}^h$, where so-called *primal continuity constraints* hold throughout the iterations. Those primal constraints are defined such that:

1. each subdomain problem becomes invertible (which is not the case when working in $\mathcal{W}^h$),

2. good convergence bound can be obtained for the DD solver.

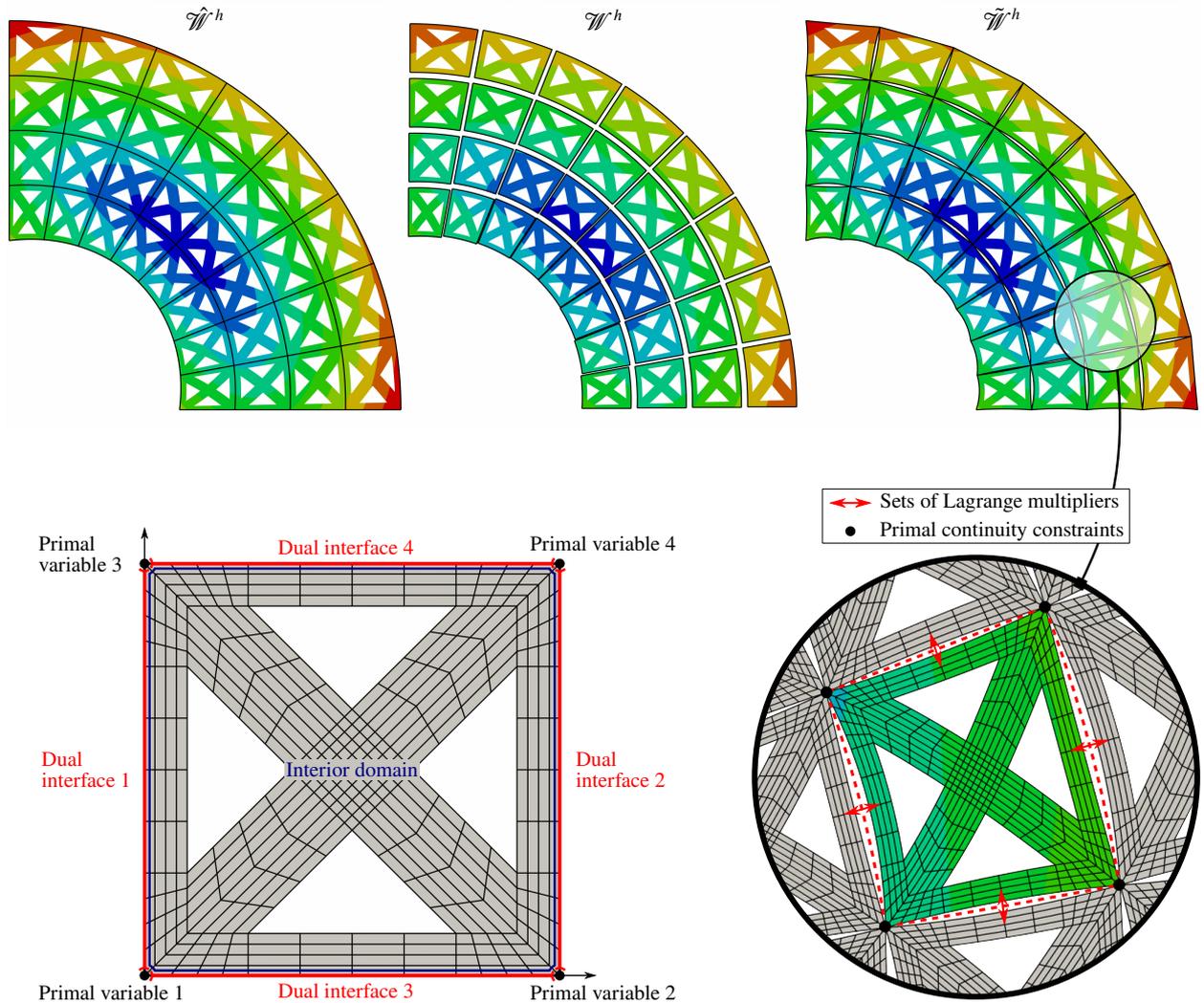

Figure 3: Illustration of the different solution spaces for the DD algorithms when applied to lattice structures. In our case of a FETI-DP based algorithm, the solution is sought in the intermediate space $\widetilde{\mathscr{W}}^h$, where the continuity of the local displacements is enforced at the corners of the subdomains.

We denote this space as $\widetilde{\mathscr{W}}^h$, for which we have $\widehat{\mathscr{W}}^h \subset \widetilde{\mathscr{W}}^h \subset \mathscr{W}^h$. It can be defined, to some extent, as:

$$\widetilde{\mathscr{W}}^h = \left\{ u^h \in \mathscr{W}^h : u^h \text{ verifies the primal continuity constraints across } \Gamma \right\}. \tag{10}$$

The introduced primal continuity constraints only affect a subset of interface variables (defined as *primal variables*). The primal continuity constraints depend on the spatial dimension and on the heterogeneity of the problem that is considered (see, *e.g.*, [37] for a theoretical study). In this work, we will simply consider the corners of the subdomains for the primal variables [18, 46], which appears sufficient in our numerical experiments to obtain a scalable method. An illustration of the different solution spaces is given in Figure 3. We can eventually split the interface variables into two categories, the primal and *dual variables*, and, then, the FETI-DP space $\widetilde{\mathscr{W}}^h$ can be viewed as:

$$\widetilde{\mathscr{W}}^h = \mathscr{W}_I^h \oplus \widetilde{\mathscr{W}}_\Gamma^h, \quad \text{with} \quad \widetilde{\mathscr{W}}_\Gamma^h = \mathscr{W}_D^h \oplus \widehat{\mathscr{W}}_P^h. \tag{11}$$

Hence, the interface space is composed of the continuous members of the primal space (*i.e.*, space $\widehat{\mathscr{W}}_P^h$) and the standard dual space (*i.e.*, space $\mathscr{W}_D^h$) with possibly disconnected (or discontinuous) members.

Additionally, Lagrange multipliers are introduced in order to recover the continuity over the dual space along the iterative algorithm. Thus, we enforce continuity for all dual interface variables in a weak sense. As we placed ourselves in a fully matching setting (since the same reference mesh is used for each cell, see Figures 2 and 3), we can directly measure the interface jump by considering the differences between couples of displacement Degrees Of Freedom (DOF), and assign a Lagrange multiplier to each of them to introduce a continuity constraint. Importantly, we also enforce the Dirichlet boundary conditions over the dual variables through Lagrange multipliers in order to keep local dual spaces with the same cardinality.

### 2.2.3. Matrix formulation and Saddle point problem

Let us now highlight how the imposition of the Primal-Dual coupling constraints is done in practice by employing a matrix notation. First, we partition every set of local displacement DOF into three categories: the internal, the dual, and the primal DOF denoted with subscripts $\cdot_i$, $\cdot_d$, and $\cdot_p$, respectively. This reads:

$$\mathbf{u}^{(s)} = \left(\mathbf{u}_i^{(s)}, \mathbf{u}_d^{(s)}, \mathbf{u}_p^{(s)}\right), \quad \text{with} \quad \mathbf{u}_i^{(s)} \in \mathbb{R}^{n_i}, \mathbf{u}_d^{(s)} \in \mathbb{R}^{n_d}, \mathbf{u}_p^{(s)} \in \mathbb{R}^{n_p}. \tag{12}$$

We additionally denote by $\mathbf{u}_r^{(s)} = (\mathbf{u}_i^{(s)}, \mathbf{u}_d^{(s)}) \in \mathbb{R}^{n_r}$ the non-primal DOF (also called the remaining DOF), and we introduce a restriction operator $\mathbf{T}_{dr} \in \mathbb{R}^{n_d \times n_r}$ that enables to extract the dual DOF from the remaining ones as $\mathbf{u}_d^{(s)} = \mathbf{T}_{dr}\mathbf{u}_r^{(s)}$. All the local primal unknowns are assembled globally into $\mathbf{u}_P$, and the remaining ones are simply concatenated which leads to define the total vector of unknowns $\mathbf{u} \in \mathbb{R}^{n_U}$ as:

$$\mathbf{u} = \left(\mathbf{u}_r^{(1)}, \mathbf{u}_r^{(2)}, \ldots, \mathbf{u}_r^{(N)}, \mathbf{u}_P\right). \tag{13}$$

In the following we will also denote by $n_R$ and $n_P$ the total number of remaining and primal DOF ($n_U = n_R + n_P$), respectively.

In order to enforce the continuity between the dual interfaces, we introduce a (discrete) jump operator of the following form:

$$B : \tilde{\mathcal{W}}^h \to \mathbb{R}^L; u^h \mapsto B(u^h) = \mathbf{Bu} - \mathbf{d}, \tag{14}$$

where $L$ is the total number of discrete continuity equations and Dirichlet constraints. Vector $\mathbf{d} \in \mathbb{R}^L$ has zero entries where the coupling conditions concern a subdomain-to-subdomain constraint, and integrates the Dirichlet conditions (Equation (3c)) otherwise. The coupling matrix $\mathbf{B} \in \mathbb{R}^{L \times n_U}$ is formed by local (dual) Boolean operators $\mathbf{B}_d^{(s)} \in \mathbb{R}^{L \times n_d}$ which are rectangular sparse matrices with only $0, 1, -1$ entries. More precisely, it reads:

$$\mathbf{B} = \begin{bmatrix} \mathbf{B}_R & \mathbf{0}_{LP} \end{bmatrix} \quad \text{with} \quad \mathbf{B}_R = \begin{bmatrix} \mathbf{B}_r^{(1)} & \ldots & \mathbf{B}_r^{(N)} \end{bmatrix}, \quad \text{and} \quad \mathbf{B}_r^{(s)} = \mathbf{B}_d^{(s)}\mathbf{T}_{dr}. \tag{15}$$

The interested readers may refer to [20, 32] to find details on the construction of these operators.

Finally, we define a set of pointwise Lagrange multipliers $\lambda \in \mathbb{R}^L$ to enforce the continuity constraints $\mathbf{Bu} - \mathbf{d} = \mathbf{0}$. It leads to formulate the following saddle problem:

$$\begin{bmatrix} \mathbf{K} & \mathbf{B}^\top \\ \mathbf{B} & \mathbf{0} \end{bmatrix} \begin{pmatrix} \mathbf{u} \\ \lambda \end{pmatrix} = \begin{pmatrix} \mathbf{f} \\ \mathbf{d} \end{pmatrix}, \tag{16}$$

where, due to the dual-primal domain decomposition, we have:

$$\mathbf{K} = \begin{bmatrix} \mathbf{K}_{RR} & \mathbf{K}_{RP} \\ \mathbf{K}_{PR} & \mathbf{K}_{PP} \end{bmatrix}, \quad \mathbf{f} = \begin{pmatrix} \mathbf{f}_R \\ \mathbf{f}_P \end{pmatrix}. \tag{17}$$

As already stated, the primal DOF are assembled globally within FETI-DP approaches. This is a common task in finite elements and can be done by employing a global-local assembly operator $\mathbf{A}_p^{(s)} \in \mathbb{R}^{n_P \times n_p}$ per subdomain. Therefore, the quantities in Equation (17) are obtained by:

$$\mathbf{K}_{PP} = \sum_{s=1}^{N} \mathbf{A}_p^{(s)} \mathbf{K}_{pp}^{(s)} \mathbf{A}_p^{(s)\top}, \quad \mathbf{f}_P = \sum_{s=1}^{N} \mathbf{A}_p^{(s)} \mathbf{f}_p^{(s)}, \tag{18}$$

and:

$$\mathbf{K}_{RR} = \begin{bmatrix} \mathbf{K}_{rr}^{(1)} & & \\ & \ddots & \\ & & \mathbf{K}_{rr}^{(N)} \end{bmatrix}, \quad \mathbf{K}_{RP} = \begin{bmatrix} \mathbf{K}_{rP}^{(1)} \\ \vdots \\ \mathbf{K}_{rP}^{(N)} \end{bmatrix}, \quad \text{with} \quad \mathbf{K}_{rP}^{(s)} = \mathbf{K}_{rp}^{(s)} \mathbf{A}_p^{(s)\top}, \quad \mathbf{K}_{PR} = \mathbf{K}_{RP}^{\top}. \tag{19}$$

Hence, regarding the structure of the total matrix $\mathbf{K}$, it is interesting to notice the presence of the block diagonal term $\mathbf{K}_{RR}$ made of the local problems $\mathbf{K}_{rr}^{(s)}$. In addition to being symmetric positive definite, this term is preponderant as $n_R \gg n_P$ (since the number of primal constraints is usually low). Thus, the main idea behind FETI-DP approaches (and DD methods in general), is to take advantage of this particular structure.

### 2.2.4. FETI-DP algorithm

The common starting point of FETI solvers consists in removing the displacement unknowns from system (16) in order to formulate an interface problem where only the Lagrange multipliers remain. With the FETI-DP approach, this is done by firstly condensing the remaining DOF by using the relation:

$$\mathbf{u}_R = \mathbf{K}_{RR}^{-1}\big(\mathbf{f}_R - \mathbf{K}_{RP}\mathbf{u}_P - \mathbf{B}_R^{\top}\lambda\big), \tag{20}$$

and then the primal DOF by using the relation:

$$\mathbf{u}_P = \mathbf{S}_{PP}^{-1}\big((\mathbf{f}_P - \mathbf{K}_{PR}\mathbf{K}_{RR}^{-1}\mathbf{f}_R) + [\mathbf{K}_{PR}\mathbf{K}_{RR}^{-1}\mathbf{B}_R^{\top}]\lambda\big), \tag{21}$$

where $\mathbf{S}_{PP} \in \mathbb{R}^{n_P \times n_P}$ denotes the primal Schur complement. It consists in the assembly of all the local primal Schur complements $\mathbf{S}_{pp}^{(s)}$:

$$\mathbf{S}_{PP} = \sum_{s=1}^{N} \mathbf{A}_p^{(s)} \mathbf{S}_{pp}^{(s)} \mathbf{A}_p^{(s)\top}, \quad \text{with} \quad \mathbf{S}_{pp}^{(s)} = \mathbf{K}_{pp}^{(s)} - \mathbf{K}_{pr}^{(s)} \mathbf{K}_{rr}^{(s)-1} \mathbf{K}_{rp}^{(s)}, \tag{22}$$

and can be written in a more compact way as: $\mathbf{S}_{PP} = \mathbf{K}_{PP} - \mathbf{K}_{PR}\mathbf{K}_{RR}^{-1}\mathbf{K}_{RP}$.

Using equations (20) and (21) enables to transform the continuity constraint between the dual interfaces (*i.e.*, the second equation in the saddle point problem in Equation (16)) into the following interface problem:

$$\mathbf{Bu} = \mathbf{d} \quad \Leftrightarrow \quad \mathbf{F}\lambda = \bar{\mathbf{d}}, \tag{23}$$

with:

$$\mathbf{F} = \mathbf{B}_R \mathbf{K}_{RR}^{-1} \mathbf{B}_R^{\top} + \mathbf{B}_R \big[\mathbf{K}_{RR}^{-1} \mathbf{K}_{RP} \mathbf{S}_{PP}^{-1} \mathbf{K}_{PR} \mathbf{K}_{RR}^{-1}\big] \mathbf{B}_R^{\top}, \tag{24}$$

$$\bar{\mathbf{d}} = -\mathbf{d} + \mathbf{B}_R \mathbf{K}_{RR}^{-1} \mathbf{f}_R - \mathbf{B}_R \mathbf{K}_{RR}^{-1} \mathbf{K}_{RP} \mathbf{S}_{PP}^{-1} \big(\mathbf{f}_P - \mathbf{K}_{PR} \mathbf{K}_{RR}^{-1} \mathbf{f}_R\big). \tag{25}$$

$\mathbf{F}$, which is called the global dual Schur complement, may also be expressed in a more compact form as:

$$\mathbf{F} = \mathbf{B}\mathbf{K}^{-1}\mathbf{B}^{\top}. \tag{26}$$

The global dual Schur complement is made of two terms as shown in Equation (24). The first term collects the local Neumann schur complements $\mathbf{F}_{dd}^{(s)}$:

$$\mathbf{B}_R \mathbf{K}_{RR}^{-1} \mathbf{B}_R^{\top} = \sum_{s=1}^{N} \mathbf{B}_d^{(s)} \mathbf{F}_{dd}^{(s)} \mathbf{B}_d^{(s)\top}, \quad \text{with} \quad \mathbf{F}_{dd}^{(s)} = \mathbf{T}_{dr} \mathbf{K}_{rr}^{(s)-1} \mathbf{T}_{dr}^{\top}, \tag{27}$$

and the second term is the coarse problem correction associated to the primal variables. Let us further detail key quantities in Equations (24) and (25) that involve the local dual and primal solutions denoted $\mathbf{U}_{rd}^{(s)}$ and $\mathbf{U}_{rp}^{(s)}$:

$$\mathbf{B}_R \mathbf{K}_{RR}^{-1} = \begin{bmatrix} \mathbf{B}_d^{(1)} \mathbf{U}_{rd}^{(1)\top} & \cdots & \mathbf{B}_d^{(N)} \mathbf{U}_{rd}^{(N)\top} \end{bmatrix}, \quad \text{with} \quad \mathbf{U}_{rd}^{(s)} = \mathbf{K}_{rr}^{(s)-1} \mathbf{T}_{dr}^{\top}, \tag{28}$$

$$\mathbf{K}_{RR}^{-1} \mathbf{K}_{RP} = \begin{bmatrix} \mathbf{U}_{rp}^{(1)} \mathbf{A}_p^{(1)\top} \\ \vdots \\ \mathbf{U}_{rp}^{(N)} \mathbf{A}_p^{(N)\top} \end{bmatrix}, \quad \text{with} \quad \mathbf{U}_{rp}^{(s)} = \mathbf{K}_{rr}^{(s)-1} \mathbf{K}_{rp}^{(s)}. \tag{29}$$

Secondly, the main idea is to solve the interface problem (23) with an iterative solver, as for instance the preconditioned conjugate gradient since $\mathbf{F}$ is symmetric positive definite. This way, only matrix-vector products need to be performed. This can be done without building the dual Schur complement $\mathbf{F}$. Indeed, in practice only the local primal solutions, the local primal Schur complements, and the global Schur complement are explicitly built. Furthermore, the application of $\mathbf{F}$ to a vector involves the solution of local independent problems, which can be done in parallel.

The last thing to do thus consists in defining a preconditioner for system (23). Usually, the so-called Dirichlet preconditioner is considered in linear elasticity [56]. It is built as follows:

$$\mathbf{M}_\mathrm{D}^{-1} = \sum_{s=1}^{N} \mathbf{D}^{(s)} \mathbf{B}_d^{(s)} \mathbf{S}_{dd}^{(s)} \mathbf{B}_d^{(s)\top} \mathbf{D}^{(s)}, \quad \text{with} \quad \mathbf{S}_{dd}^{(s)} = \mathbf{K}_{dd}^{(s)} - \mathbf{K}_{di}^{(s)} \mathbf{K}_{ii}^{(s)-1} \mathbf{K}_{id}^{(s)}, \qquad (30)$$

where $\mathbf{D}^{(s)} \in \mathbb{R}^{L \times L}$ are diagonal weight matrices that verify $\sum_{s=1}^{N} \mathbf{B}_d^{(s)} \mathbf{B}_d^{(s)\top} \mathbf{D}^{(s)} = \mathbf{I}_L$ and can be defined in order to better balance the interface correction in case of material heterogeneities and redundancies in the compatibility constraints, as explained in, *e.g.*, [20]. More elaborated strategies exist to handle very high heterogeneities, as for instance, the use of Generalized Eigenproblems in the Overlaps [62] or multipreconditioned Krylov solvers [9, 61]. Once the interface problem (23) is solved, one can recover the displacement solution $\mathbf{u}$ by using successively Equations (21) and (20).

Overall, returning to the non-condensed form (16) and making use of (26), it is possible to rewrite this domain decomposition solution in a matrix from as follows:

$$\begin{pmatrix} \mathbf{u} \\ \boldsymbol{\lambda} \end{pmatrix} = \underbrace{\begin{bmatrix} \mathbf{I} & -\mathbf{K}^{-1} \mathbf{B}^\top \\ & \mathbf{I} \end{bmatrix}}_{\text{Recover full sol.}} \underbrace{\begin{bmatrix} \mathbf{K}^{-1} & \\ & -\mathbf{F}^{-\mathtt{pcg}} \end{bmatrix}}_{\text{Interface pb.}} \underbrace{\begin{bmatrix} \mathbf{I} & \\ -\mathbf{B}\mathbf{K}^{-1} & \mathbf{I} \end{bmatrix}}_{\text{Condense}} \begin{pmatrix} \mathbf{f} \\ \mathbf{d} \end{pmatrix}, \qquad (31)$$

where $\mathbf{F}^{-\mathtt{pcg}}$ stands for the application of $\mathbf{F}^{-1}$ via a preconditioned conjugate gradient (*i.e.*, solving the interface problem).

### 2.3. Block preconditioners saddle point problems

Applying the standard algorithm explained in Section 2.2.4 requires to solve *exactly* (*i.e.*, close to machine precision) the local problems (within $\mathbf{K}_{RR}^{-1}$) and the coarse problem ($\mathbf{S}_{PP}^{-1}$) at every iteration, see Equations (24) and (25). In practice, it is commonly done by factorizing all these operators *a priori*, and then performing backward-forward substitutions during the iterative algorithm. It is worth nothing that $\mathbf{K}_{RR}^{-1}$ involves the solution of local independent problems; therefore, it can be done in parallel. However, this still raises issues, for instance, in the case of large coarse problem (or large local problem) which leads to high memory usage (and bad parallel efficiency). In the specific case of lattice structures, one may seek to exploit the similarities between the local problems. If exact local solutions are mandatory, then it seems not straightforward to take advantage of these similarities (unless all cells are identical).

Therefore, we will rather consider an alternative approach that consists in working directly with the saddle point problem (16), and to use a suitable iterative solver (GMRES for example) together with a block preconditioner. It is possible to build block preconditioners dedicated to the saddle point problem (16) that uses the same ingredients as those involved in standard FETI approaches (local and coarse solutions). However, since the treatment now appears in the preconditioner, it is possible to approximate the local problems (and also the coarse problem), solving them *inexactly*. This leads to the so-called inexact FETI-DP approaches [7, 36, 55].

The starting point to build block preconditioners for saddle point problems of type (16) consists in performing the following block $\mathrm{LDL}^\top$ factorization:

$$\begin{bmatrix} \mathbf{K} & \mathbf{B}^\top \\ \mathbf{B} & \mathbf{0} \end{bmatrix} = \begin{bmatrix} \mathbf{I} & \\ \mathbf{B}\mathbf{K}^{-1} & \mathbf{I} \end{bmatrix} \begin{bmatrix} \mathbf{K} & \\ & -\mathbf{F} \end{bmatrix} \begin{bmatrix} \mathbf{I} & \mathbf{K}^{-1} \mathbf{B}^\top \\ & \mathbf{I} \end{bmatrix}, \qquad (32)$$

which is possible as $\mathbf{K}$ is a symmetric positive definite matrix and $\mathbf{B}$ has full rank. Then, standard block preconditioners consist in using different sub-preconditioners, *i.e.*, different approximations of each blocks, see for instance Notay [48].

Indeed, block preconditioners are of the form:

$$\mathring{\mathbf{A}}^{-1} = \begin{bmatrix} \mathbf{I} & -\mathring{\mathbf{U}}^\top \\ & \mathbf{I} \end{bmatrix} \begin{bmatrix} \mathring{\mathbf{K}}^{-1} & \\ & -\mathring{\mathbf{F}}^{-1} \end{bmatrix} \begin{bmatrix} \mathbf{I} & \\ -\mathring{\mathbf{V}} & \mathbf{I} \end{bmatrix}, \tag{33}$$

where $\mathring{\mathbf{K}}^{-1}$ is a preconditioner for $\mathbf{K}$, $\mathring{\mathbf{F}}^{-1}$ is a preconditioner for $\mathbf{F}$, and $\mathring{\mathbf{U}}$ and $\mathring{\mathbf{V}}$ are approximations of $\mathbf{B}\mathbf{K}^{-1}$. Several choices for these sub-preconditioners can be considered. In the context of FETI-based approaches, we could select:

- Standard FETI-DP (see Equation (31)): $\mathring{\mathbf{K}}^{-1} = \mathbf{K}^{-1}$, $\mathring{\mathbf{F}}^{-1} = \mathbf{F}^{-\texttt{pcg}}$, $\mathring{\mathbf{U}} = \mathring{\mathbf{V}} = \mathbf{B}\mathbf{K}^{-1}$,

- Block triangular Inexact FETI-DP from Klawonn and Rheinbach [36]:
  $\mathring{\mathbf{K}}^{-1} = \mathbf{K}^{-\texttt{AMG}}$, $\mathring{\mathbf{F}}^{-1} = \mathbf{M}_{\text{D}}^{-\texttt{AMG}}$, $\mathring{\mathbf{U}} = \mathbf{0}$, $\mathring{\mathbf{V}} = \mathbf{B}\mathbf{K}^{-1}$,
  where $(\cdot)^{-\texttt{AMG}}$ means that the solution is obtained with an Algebraic MultiGrid method,

- Block diagonal Inexact IETI-DP based on Bosy et al. [7]:
  $\mathring{\mathbf{K}}^{-1} = \mathbf{K}^{-\texttt{FD}}$, $\mathring{\mathbf{F}}^{-1} = \mathbf{M}_{\text{D}}^{-\texttt{FD}}$, $\mathring{\mathbf{U}} = \mathbf{0}$, $\mathring{\mathbf{V}} = \mathbf{0}$,
  where $(\cdot)^{-\texttt{FD}}$ means that the inverses are approximated via the Fast Diagonalization method.

In this work, we seek to introduce sub-preconditioners in Equation (33) that exploit the similarities between all local subproblems, *i.e.*, the fact that:

$$\mathbf{K}^{(1)} \simeq \mathbf{K}^{(2)} \simeq \ldots \simeq \mathbf{K}^{(N)}. \tag{34}$$

The idea is therefore to use only a small set of local stiffness matrices to efficiently approximate all the local solutions, which is the object of next section.

## 3. Benefiting from the repetitiveness of lattices

### 3.1. Reduced Order Modeling

#### 3.1.1. Reduced Bases and Principal subdomains

Let us now introduce the matrix manifold comprising all stiffness matrices associated to every sub-problem of our total DD problem:

$$\mathscr{M}_{\mathbf{K}} = \left\{ \mathbf{K}^{(1)}, \ldots, \mathbf{K}^{(N)} \right\}. \tag{35}$$

Due to the similarities between these matrices, we suspect the manifold $\mathscr{M}_{\mathbf{K}}$ to be of low dimension; *i.e.*, it can be well represented by the span of a low number $N_{\text{rb}} \ll N$ of reference matrices, denoted as $\mathbf{K}^{[k]}$, with $k = 1, \ldots, N_{\text{rb}}$. Thus, we assume that there is a Reduced Basis (RB) space:

$$\mathbb{K} = \text{span}\left\{ \mathbf{K}^{[1]}, \ldots, \mathbf{K}^{[N_{\text{rb}}]} \right\}, \tag{36}$$

such that:

$$\forall \mathbf{K}^{(s)} \in \mathscr{M}_{\mathbf{K}}, \quad \exists \mathring{\mathbf{K}}^{(s)} \in \mathbb{K} \ s.t. \ \|\mathbf{K}^{(s)} - \mathring{\mathbf{K}}^{(s)}\|_v \leq \text{tol}, \tag{37}$$

where $\|\cdot\|_v$ is an appropriate norm. In Section 4.2, we will explain in detail how to build efficiently the reduced space $\mathbb{K}$. At this stage, let us just mention that we will make use of a greedy approach which generates an interpolatory basis; that is, the reference matrices $\mathbf{K}^{[k]}$ correspond to local matrices $\mathbf{K}^{(s)}$ of appropriately selected subdomains denoted $\sigma_k, k = 1, \ldots, N_{\text{rb}}$. In other words:

$$\mathring{\mathbf{K}}^{(s)} = \sum_{k=1}^{N_{\text{rb}}} \alpha_k^{(s)} \mathbf{K}^{(\sigma_k)}, \tag{38}$$

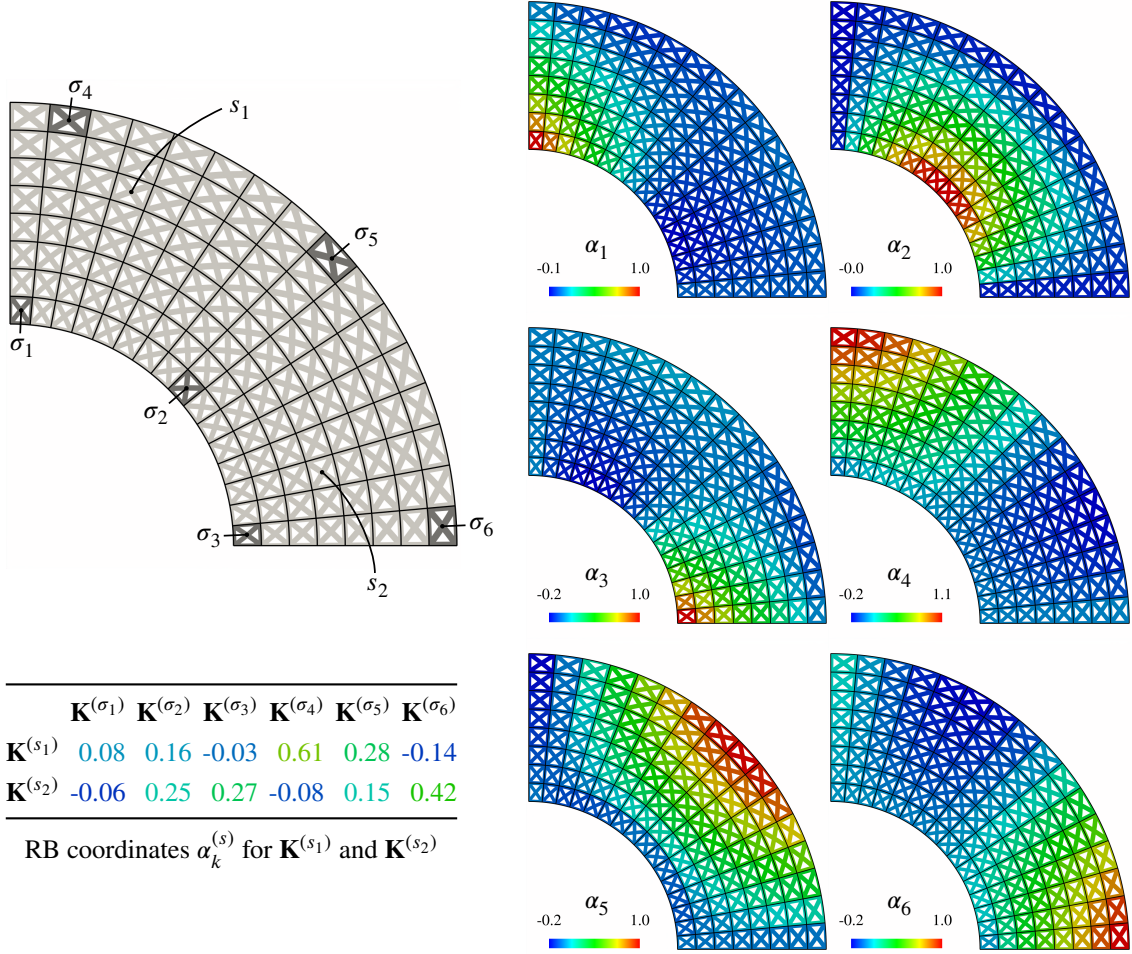

Figure 4: Example of results obtained with the reduced basis approach. 6 principal subdomains (marked as $\sigma_k, k = 1, \ldots, 6$) and associated local stiffness matrices are selected (see top-left). Then, all the local stiffness matrices can be expressed as a linear combination of the principal ones. The exact values of the associated RB coordinates are given for subdomains $s_1$ and $s_2$ in the table at the bottom-left. For completeness, plots of the value of each RB coordinate for all subdomains are realized on the right.

where the $\alpha_k^{(s)} \in \mathbb{R}$ are coefficients which depend on the subdomain $s$. Figure 4 illustrates the result of such a process, and, once again, Section 4.2 will detail the construction of these approximations, including the selection of the principal subdomains.

Similarities among subdomains are also reflected in all the different local quantities involved in the FETI-DP method (see again Section 2.2.4), *i.e.*, in:

- Primal solution: $\mathbf{U}_{rp}^{(s)} = \mathbf{K}_{rr}^{(s)-1} \mathbf{K}_{rp}^{(s)}$,
- Primal Schur: $\mathbf{S}_{pp}^{(s)} = \mathbf{K}_{pp}^{(s)} - \mathbf{K}_{pr}^{(s)} \mathbf{K}_{rr}^{(s)-1} \mathbf{K}_{rp}^{(s)}$,
- Dual solution: $\mathbf{U}_{rd}^{(s)} = \mathbf{K}_{rr}^{(s)-1} \mathbf{T}_{rd}$,
- Neumann Schur: $\mathbf{F}_{dd}^{(s)} = \mathbf{T}_{dr} \mathbf{K}_{rr}^{(s)-1} \mathbf{T}_{rd}$,
- Dirichlet Schur: $\mathbf{S}_{dd}^{(s)} = \mathbf{K}_{dd}^{(s)} - \mathbf{K}_{di}^{(s)} \mathbf{K}_{ii}^{(s)-1} \mathbf{K}_{id}^{(s)} = \mathbf{F}_{dd}^{(s)-1}$.

Based on the reduced basis of the local stiffness matrices and the selection of principal subdomains, we also consider reduced bases of all these local quantities by simply replacing $\mathbf{K}^{(s)}$ by $\mathbf{K}^{(\sigma_k)}$ in the above equations, which leads to:

$$\mathbb{U}_{rp} = \text{span}\left\{\mathbf{U}_{rp}^{(\sigma_1)}, \ldots, \mathbf{U}_{rp}^{(\sigma_N)}\right\}, \quad \mathbb{S}_{pp} = \text{span}\left\{\mathbf{S}_{pp}^{(\sigma_1)}, \ldots, \mathbf{S}_{pp}^{(\sigma_N)}\right\}, \tag{39a}$$

$$\mathbb{U}_{rd} = \text{span}\left\{\mathbf{U}_{rd}^{(\sigma_1)}, \ldots, \mathbf{U}_{rd}^{(\sigma_N)}\right\}, \quad \mathbb{F}_{dd} = \text{span}\left\{\mathbf{F}_{dd}^{(\sigma_1)}, \ldots, \mathbf{F}_{dd}^{(\sigma_N)}\right\}, \quad \mathbb{S}_{dd} = \text{span}\left\{\mathbf{S}_{dd}^{(\sigma_1)}, \ldots, \mathbf{S}_{dd}^{(\sigma_N)}\right\}. \tag{39b}$$

Then, the idea is to use these RB in the inexact FETI-DP preconditioner by approximating all the local operators.

**Remark 2.** *All the reduced bases for the DD operators have the same size which is inherited from the reduced basis of the stiffness matrices* (36). *Indeed, we do not place ourselves in a standard ROM configuration where the reduced bases are built in an offline phase. We cannot pay the price of generating snapshots of all local DD operators as it would take more time than simply solving the FE system. Instead, we build these operators for the principal cells involved in $\mathbb{K}$ and use them to generate the other reduced bases (see again Equation* (39)). *Therefore, we will call upon a (global) iterative process to recover the exact solution (instead of in one go as for the standard FETI-DP approach, see Equation* (31)).

### 3.1.2. Approximating the local operators

Let us invoke the main ideas behind ROM [25, 54] in order to benefit from the RB spaces. We project each subdomain local operator into their corresponding reduced bases, or, in other words, we aim at finding the coefficients of the local operators in the reduced bases.

To this purpose, we solve two minimization problems per subdomain, which enable to approximate the primal and dual solutions denoted previously $\mathbf{U}_{rp}^{(s)}$ and $\mathbf{U}_{rd}^{(s)}$, respectively:

$$\mathring{\mathbf{U}}_{rp}^{(s)} = \underset{\mathbf{U}_{rp} \in \mathbb{U}_{rp}}{\text{argmin}} \left\|\mathbf{U}_{rp}^{(s)} - \mathbf{U}_{rp}\right\|_{E^{(s)}}^{2} \quad \text{and} \quad \mathring{\mathbf{U}}_{rd}^{(s)} = \underset{\mathbf{U}_{rd} \in \mathbb{U}_{rd}}{\text{argmin}} \left\|\mathbf{U}_{rd}^{(s)} - \mathbf{U}_{rd}\right\|_{E^{(s)}}^{2}, \tag{40}$$

where $\|\cdot\|_{E^{(s)}}^2$ denotes an energy norm which for a vector $\mathbf{v}_r \in \mathbb{R}^r$ and a column matrix $\mathbf{V}_{rl} \in \mathbb{R}^{r \times l}$ and read as, respectively:

$$\|\mathbf{v}_r\|_{E^{(s)}}^2 = \frac{1}{2}\mathbf{v}_r^\top \mathbf{K}_{rr}^{(s)} \mathbf{v}_r, \quad \text{and} \quad \|\mathbf{V}_{rl}\|_{E^{(s)}}^2 = \frac{1}{2}\text{tr}\left(\mathbf{V}_{rl}^\top \mathbf{K}_{rr}^{(s)} \mathbf{V}_{rl}\right).$$

In case of column matrices, the energy norm can also be seen as the sum of the energy norms of each column vectors. With such a choice, solving Equation (40) is equivalent to solve small $N_{rb}$-by-$N_{rb}$ linear systems:

$$\mathbf{A}^{(s)} \boldsymbol{\pi}^{(s)} = \mathbf{b}^{(s)} \quad \text{and} \quad \mathbf{M}^{(s)} \boldsymbol{\delta}^{(s)} = \mathbf{d}, \tag{41}$$

where $\mathbf{A}^{(s)}, \mathbf{M}^{(s)} \in \mathbb{R}^{N_{rb} \times N_{rb}}$ and $\mathbf{b}^{(s)}, \mathbf{d} \in \mathbb{R}^{N_{rb}}$ are given by:

$$\mathbf{A}^{(s)}[i,j] = \text{tr}\left(\mathbf{U}_{rp}^{(\sigma_i)\top} \mathbf{K}_{rr}^{(s)} \mathbf{U}_{rp}^{(\sigma_j)}\right), \quad \mathbf{b}^{(s)}[i] = \text{tr}\left(\mathbf{U}_{rp}^{(\sigma_i)\top} \mathbf{K}_{rp}^{(s)}\right), \tag{42a}$$

$$\mathbf{M}^{(s)}[i,j] = \text{tr}\left(\mathbf{U}_{rd}^{(\sigma_i)\top} \mathbf{K}_{rr}^{(s)} \mathbf{U}_{rd}^{(\sigma_j)}\right), \quad \mathbf{d}[i] = \text{tr}\left(\mathbf{U}_{rd}^{(\sigma_i)\top} \mathbf{T}_{rd}\right). \tag{42b}$$

We used the syntax $[i,j]$ and $[i]$ for the index notation in order to distinguish it from the size/DOF related notation needed for the description of the FETI algorithm.

Once the coefficients $\boldsymbol{\alpha}^{(s)}$ (defined in Equation (38)), $\boldsymbol{\pi}^{(s)}$, and $\boldsymbol{\delta}^{(s)}$ are known, we can obtain the ROM of all the local DD operators as:

$$\mathring{\mathbf{U}}_{rp}^{(s)} = \mathbb{U}_{rp} \boldsymbol{\pi}^{(s)}, \quad \mathring{\mathbf{S}}_{pp}^{(s)} = \text{sym}\left(\mathring{\mathbf{K}}_{pp}^{(s)} - \mathring{\mathbf{K}}_{rp}^{(s)\top} \mathring{\mathbf{U}}_{rp}^{(s)}\right), \quad \mathring{\mathbf{U}}_{rd}^{(s)} = \mathbb{U}_{rd} \boldsymbol{\delta}^{(s)}, \quad \mathring{\mathbf{F}}_{dd}^{(s)} = \mathbb{F}_{dd} \boldsymbol{\delta}^{(s)}, \quad \mathring{\mathbf{S}}_{dd}^{(s)} = \mathbb{S}_{dd} \boldsymbol{\alpha}^{(s)}. \tag{43}$$

In the above Equation, we made an abuse of notation for the sake of clarity: For instance, $\mathbb{U}_{rp} \boldsymbol{\pi}^{(s)}$ should be understood as the linear combination of the basis elements spanning the space $\mathbb{U}_{rp}$ with the scalars collected in vector $\boldsymbol{\pi}^{(s)}$. And the same applies to $\mathbb{U}_{rd} \boldsymbol{\delta}^{(s)}$, $\mathbb{F}_{dd} \boldsymbol{\delta}^{(s)}$, and $\mathbb{S}_{dd} \boldsymbol{\alpha}^{(s)}$.

**Algorithm 1** Steps when applying the iFETI-DP preconditioner to an input vector.

Input vector:

In: $\begin{pmatrix} \mathbf{v} \\ \mathbf{w} \end{pmatrix}$, with $\mathbf{v} \in \mathbb{R}^{n_U}$, $\mathbf{w} \in \mathbb{R}^L$.

---

| | | | |
|---|---|---|---|
| $\bar{\mathbf{w}}$ | $\leftarrow$ | $\mathbf{w} - \mathring{\mathbf{U}}\mathbf{v}$ | *Condense the right-hand side (see Section 3.2.4)* |
| $\mathbf{y}$ | $\leftarrow$ | $-\mathring{\mathbf{F}}^{-\text{pcg}}\bar{\mathbf{w}}$ | *Solve the interface problem (see Section 3.2.3)* |
| $\mathbf{x}$ | $\leftarrow$ | $\mathring{\mathbf{K}}^{-\text{ROM}}\mathbf{v} - \mathring{\mathbf{U}}^\top \mathbf{y}$ | *Recover the full solution (see Section 3.2.2)* |

Ouput vector:

Out: $\begin{pmatrix} \mathbf{x} \\ \mathbf{y} \end{pmatrix} = \mathring{\mathbf{A}}_{\text{iFETIdp}}^{-1} \begin{pmatrix} \mathbf{v} \\ \mathbf{w} \end{pmatrix}$, with $\mathbf{x} \in \mathbb{R}^{n_U}$, $\mathbf{y} \in \mathbb{R}^L$.

---

**Remark 3.** *The local Dirichlet Schur complements $\mathbf{S}_{dd}^{(s)}$ are approximated directly through the ROM coefficients $\boldsymbol{\alpha}^{(s)}$ of the stiffness matrices. We could improve this approximation by employing an expression similar to the one of the local primal Schur complements. However, this would add significant extra computational costs (projection of $\mathbf{K}_{ii}^{(s)-1}\mathbf{K}_{id}^{(s)}$ using a similar approach as in (40)) for little improvements as these operators are involved in the preconditioning of the interface problem only (see the expression of the Dirichlet preconditioner in Equation (30)). Our numerical experiments confirm that it is a sufficient choice.*

**Remark 4.** *Contrary to Remark 3, we could have chosen to approximate the local primal Schur complements as $\mathring{\mathbf{S}}_{pp}^{(s)} = \mathbb{S}_{pp}\boldsymbol{\alpha}^{(s)}$ as done for the local Dirichlet Schur complements. However, since we already paid the price to approximate the primal solutions properly by determining the coefficients $\boldsymbol{\pi}^{(s)}$, it appears natural to use these best approximations.*

## 3.2. ROM-based inexact FETI-DP

We finally combine the FETI-DP approach together with the previously introduced ROM ingredients into the following block preconditioner:

$$\mathring{\mathbf{A}}_{\text{iFETIdp}}^{-1} = \begin{bmatrix} \mathbf{I} & -\mathring{\mathbf{U}}^\top \\ & \mathbf{I} \end{bmatrix} \begin{bmatrix} \mathring{\mathbf{K}}^{-\text{ROM}} & \\ & -\mathring{\mathbf{F}}^{-\text{pcg}} \end{bmatrix} \begin{bmatrix} \mathbf{I} & \\ -\mathring{\mathbf{U}} & \mathbf{I} \end{bmatrix}, \tag{44}$$

where $(\cdot)^{-\text{ROM}}$ indicates that the solution is obtained in a ROM-based, and $(\cdot)^{-\text{pcg}}$ indicates that the solution is obtained via a preconditioned conjugate gradient (as it was the case for the standard FETI-DP approach, see again Equation (31)). Algorithm 1 shows the main steps involved during the application of this inexact FETI-DP preconditioner to an input vectors, and the following sections detail the different sub-blocks.

### 3.2.1. Coarse problem

The global primal Schur complement $\mathbf{S}_{PP}$ (which appears in both $\mathring{\mathbf{K}}^{-\text{ROM}}$ and $\mathring{\mathbf{F}}^{-\text{pcg}}$) is approximated using the ROM of the local primal Schur complements $\mathring{\mathbf{S}}_{pp}^{(s)}$ as defined in Equation (43). In the offline phase, once assembled, we factorize it so that only forward-backward substitutions are performed during the iterations of the solver (*i.e.*, whenever one needs to solve the coarse problem). The factorization is done "exactly" via a Cholesky decomposition as for the standard FETI-DP approach, but we could eventually make use of approximative solvers for the coarse problem as we are in the context of inexact FETI methods. For instance, one would go for approximative solvers to maintain parallel scalabilty for large number of subdomains, as motivated in the initial works on inexact FETI-DP, see Klawonn and Rheinbach [36]. This is for sure interesting for lattice structures as they have numerous subdomains, but this has not been investigated in the scope of this work.

**Algorithm 2** Steps when applying the sub-preconditioner $\mathring{\mathbf{K}}^{-\text{ROM}}$ to an input vector.

---

Input vector:

In: $\mathbf{v} = \begin{pmatrix} \mathbf{v}_R \\ \mathbf{v}_P \end{pmatrix}$, with $\mathbf{v}_R \in \mathbb{R}^{n_R}, \mathbf{v}_P \in \mathbb{R}^{n_P}$.

---

| | | |
|---|---|---|
| $\bar{\mathbf{v}}_P$ | $\leftarrow$ $\mathbf{v}_P - \mathring{\mathbf{U}}_{RP}^\top \mathbf{v}_R$ | *Form the right-hand side of the coarse problem* |
| $\mathbf{x}_P$ | $\leftarrow$ $\mathring{\mathbf{S}}_{PP}^{-1} \bar{\mathbf{v}}_P$ | *Solve the coarse problem (see Section 3.2.1)* |
| $\mathbf{x}_R$ | $\leftarrow$ $\mathring{\mathbf{K}}_{RR}^{-\text{ROM}} \mathbf{v}_R$ | *Approximate the local problems (see Equation (46))* |
| $\mathbf{x}_R$ | $\leftarrow$ $\mathbf{x}_R - \mathring{\mathbf{U}}_{RP} \mathbf{x}_P$ | *Recover the full solution (see Section 3.2.2)* |

Ouput vector:

Out: $\mathbf{x} = \begin{pmatrix} \mathbf{x}_R \\ \mathbf{x}_P \end{pmatrix}$, with $\mathbf{x}_R \in \mathbb{R}^{n_R}, \mathbf{x}_P \in \mathbb{R}^{n_P}$.

---

### 3.2.2. Inexact local solves

The sub-preconditioner $\mathring{\mathbf{K}}^{-\text{ROM}}$ in Equation (44) reads as:

$$\mathring{\mathbf{K}}^{-\text{ROM}} = \begin{bmatrix} \mathbf{I}_R & -\mathring{\mathbf{U}}_{RP} \\ & \mathbf{I}_P \end{bmatrix} \begin{bmatrix} \mathring{\mathbf{K}}_{RR}^{-\text{ROM}} & \\ & \mathring{\mathbf{S}}_{PP}^{-1} \end{bmatrix} \begin{bmatrix} \mathbf{I}_R & \\ -\mathring{\mathbf{U}}_{PR} & \mathbf{I}_P \end{bmatrix}. \quad (45)$$

This block form originates from a block $\text{LDL}^\top$ factorization of Equation (17), as done for the complete saddle point problem (see Equation (32)). One could rewrite Equations (20) and (21) in this block format. The novelty does not lie here but in the sub-blocks.

The coarse problem is handled as described previously (see Section 3.2.1). On the other hand, the column matrix $\mathring{\mathbf{U}}_{PR}$ assembles the ROM of the local primal solutions $\mathring{\mathbf{U}}_{rp}^{(s)}$, and thus is built quickly (once the ROM is performed). Finally, operator $\mathring{\mathbf{K}}_{RR}^{-\text{ROM}}$ stands for the ROM-based solution of the local problems (recall the definition of $\mathbf{K}_{RR}$ in Equation (19)). Indeed, a matrix-vector product of type $\mathbf{x}_R = \mathring{\mathbf{K}}_{RR}^{-\text{ROM}} \mathbf{v}_R$ involves the solution of $N$ independent local problems which, for each subdomain $s$, is done similarly as the approximation of the local DD operators described in Section 3.1:

$$\mathbf{R}^{(s)} = \left[ \mathbf{K}_{rr}^{(\sigma_1)-1} \mathbf{v}_r^{(s)}, \ldots, \mathbf{K}_{rr}^{(\sigma_{N_{\mathrm{rb}}})-1} \mathbf{v}_r^{(s)} \right] \in \mathbb{R}^{r \times N_{\mathrm{rb}}}, \quad (46a)$$

$$\mathbf{x}_r^{(s)} = \underset{\mathbf{x}_r \in \mathrm{span}(\mathbf{R}^{(s)})}{\arg\min} \left\| \mathbf{K}_{rr}^{(s)-1} \mathbf{v}_r^{(s)} - \mathbf{x}_r \right\|_{E^{(s)}} = \mathbf{R}^{(s)} \left[ \mathbf{R}^{(s)\top} \mathbf{K}_{rr}^{(s)} \mathbf{R}^{(s)} \right]^{-1} \mathbf{R}^{(s)\top} \mathbf{v}_r^{(s)}. \quad (46b)$$

The first step given by Equation (46a) consists in generating a reduced basis for the local solution. This is done by using the $N_{\mathrm{rb}}$ principal local problems. Those principal stiffness matrices are previously factorized during the offline phase, and thus, the generation of the reduced bases $\mathbf{R}^{(s)}$ for the local solutions is done via forward-backward substitutions only. The second step given by Equation (46b) consists in approximating the solution within the reduced space, as already discussed in Section 3.1.2.

Algorithm 2 summarizes the main steps involved in the product between the sub-preconditioner $\mathring{\mathbf{K}}^{-\text{ROM}}$ and an input vector, and how it can be implemented in practice.

**Remark 5.** *One application of $\mathring{\mathbf{K}}^{-\text{ROM}}$ is quite expensive but presents several advantages. The main one is that only the principal local problems are involved: There is no need to factorize the N local problems as needed with standard FETI-DP (or if one need to solve the local problems exactly at every iteration) but only the $N_{\mathrm{rb}}$ principal local problems. This not only means less computational time in the preprocessing step due to the fewer factorizations, but also leads to less computer memory usage which is one of the main bottlenecks when dealing with the high-fidelity, fine-scale simulation of lattice structures.*

**Remark 6.** *Let us mention the existence of multipreconditioned iterative solvers [10, 21, 61] where the idea of generating multiple search directions per iteration is performed similarly to the process depicted in Equation (46). Indeed, one can consider using a multipreconditioned conjugate gradient solver [10, 61] to solve the local problems up to a predefined tolerance with the $\mathbf{K}^{(\sigma_k)-1}$ as the multiple preconditioners. This could lead to a $\mathring{\mathbf{K}}^{-\text{MPCG}}$ version of the subpreconditioner.*

### 3.2.3. Fast interface solves

Then, we need to define the approximation of the global dual/Neumann Schur complement used in the sub-preconditioner (*i.e.*, $\mathring{\mathbf{F}}^{-1}$ in Equation (33)). Following the inexact FETI-DP solver described in [36], one may think of using a Dirichlet preconditioner:

$$\mathring{\mathbf{F}}^{-1} = \mathring{\mathbf{M}}_D^{-\text{ROM}}, \quad (47)$$

where $\mathring{\mathbf{M}}_D^{-\text{ROM}}$ is the Dirichlet preconditioner (recall Equation (30)) but where the solutions of the local problems $\mathbf{K}_{ii}^{(s)-1}$ are approximated using a ROM-based solution similar to $\mathring{\mathbf{K}}^{-\text{ROM}}$ in Equation (46).

Although it is an interesting choice (as it reduces the number of operators to be factorized), we suggest another option in view of better approximating $\mathbf{F}^{-1}$ and thus further reducing the number of iterations of the DD solver. This is where we use the ROM of the local Schur complements $\mathring{\mathbf{F}}_{dd}^{(s)}$ and $\mathring{\mathbf{S}}_{dd}^{(s)}$ built in Section 3.1.2. By recalling the expressions of the global dual Schur complement $\mathbf{F}$ and the Dirichlet preconditioner $\mathbf{M}_D^{-1}$ defined in Equations (24) and (30), respectively, we introduce their ROM-based approximations as:

$$\mathring{\mathbf{F}} = \sum_{s=1}^{N} \mathbf{B}_d^{(s)} \mathring{\mathbf{F}}_{dd}^{(s)} \mathbf{B}_d^{(s)\top} + \mathbf{B}_R \mathring{\mathbf{U}}_{RP} \mathring{\mathbf{S}}_{PP}^{-1} \mathring{\mathbf{U}}_{RP}^\top \mathbf{B}_R^\top, \quad (48)$$

$$\mathring{\mathbf{M}}_D^{-1} = \sum_{s=1}^{N} \mathbf{D}^{(s)} \mathbf{B}_d^{(s)} \mathring{\mathbf{S}}_{dd}^{(s)} \mathbf{B}_d^{(s)\top} \mathbf{D}^{(s)}. \quad (49)$$

Beyond being an approximation, the important difference between $\mathbf{F}$ and $\mathring{\mathbf{F}}$, and $\mathbf{M}_D^{-1}$ and $\mathring{\mathbf{M}}_D^{-1}$, is the computational cost required to perform the matrix-vector product with an input vector. Indeed, as the local Schur complements of the principal subdomains are explicitly built, the application of $\mathring{\mathbf{F}}$ and $\mathring{\mathbf{M}}_D^{-1}$ to a vector only requires matrix-vector products instead of solving local problems, as requested with $\mathbf{F}$ and $\mathbf{M}_D^{-1}$. Based on this observation, we intend to use a preconditioned conjugate gradient to approximate and to solve the interfaces problem, *i.e.*, we approximate $\mathbf{F}^{-1} \approx \mathring{\mathbf{F}}^{-\text{pcg}}$, where $\mathring{\mathbf{F}}$ corresponds to (48) preconditioned with $\mathring{\mathbf{M}}_D^{-1}$, defined in Equation (49).

**Remark 7.** *The formulations of the ROM-based approximations of the global dual Schur complement and the Dirichlet preconditioner given in Equations (48) and (49) could be rewritten in terms of quantities defined over the principal subdomains. Instead of a sum of N terms, it would be a sum of $N_{\text{rb}}$ terms. The key point to remember is that those operators are cost-effective, in the sense that matrix-vector products are fast to compute.*

### 3.2.4. Fast transitions between Lagrange multipliers and displacement DOF

Finally, we approximate the rectangular matrix as $\mathring{\mathbf{U}}^\top = \mathring{\mathbf{K}}^{-\text{ROM}} \mathbf{B}^\top$ which is applied twice per call of the inexact FETI-DP preconditioner (see Equation (44)). It is involved when building the right-hand sides of the interface problems and when recovering the full solution as shown in Algorithm 1. The product between $\mathring{\mathbf{U}}$ (or $\mathring{\mathbf{U}}^\top$) and an input vector involves similar steps than the action of the sub-preconditioner $\mathring{\mathbf{K}}^{-\text{ROM}}$. The difference is that there is no need to generate reduced bases at every call. Indeed, we can use directly the ROM of the dual solutions $\mathring{\mathbf{U}}_{rd}^{(s)}$ which leads to numerically cheap matrix-vector products. Algorithm 3 summarizes the main steps involved during those products between the sub-block $\mathring{\mathbf{U}}$ and input vectors, and how it can be implemented in practice (it would be similar for $\mathring{\mathbf{U}}^\top$).

**Algorithm 3** Steps when performing matrix-vector products between $\mathring{\mathbf{U}}$ and an input vector.

Input vector:

In: $\mathbf{v} = \begin{pmatrix} \mathbf{v}_R \\ \mathbf{v}_P \end{pmatrix}$, with $\mathbf{v}_R \in \mathbb{R}^{n_R}$, $\mathbf{v}_P \in \mathbb{R}^{n_P}$.

| | | | |
|---|---|---|---|
| $\bar{\mathbf{v}}_P$ | $\leftarrow$ | $\mathbf{v}_P - \mathring{\mathbf{U}}_{RP}^\top \mathbf{v}_R$ | *Form the right-hand side of the coarse problem* |
| $\mathbf{z}_P$ | $\leftarrow$ | $\mathring{\mathbf{S}}_{PP}^{-1} \bar{\mathbf{v}}_P$ | *Solve the coarse problem (see Section 3.2.1)* |
| $\mathbf{z}$ | $\leftarrow$ | $\mathring{\mathbf{U}}_{RD}^\top \mathbf{v}_R - \mathbf{B}_R \mathring{\mathbf{U}}_{RP} \mathbf{z}_P$ | *Approximate the interface vector* |

Ouput vector:
Out: $\mathbf{z} = \mathring{\mathbf{U}} \mathbf{v}$, $\mathbf{z} \in \mathbb{R}^L$.

## 4. Application to multiscale geometric models

### 4.1. Geometric modeling via composition

We consider in this work geometric models of lattice structures based on a multiscale geometric modeling approach introduced in [16] and democratized since then (see for instance [1, 8, 28]). The idea is to model the two scales of the lattice structure separately first, and then to build the full detailed geometric model by composing the obtained micro and macro models. More specifically, within this approach, we consider that the reference microstructure is defined by a spline-based model $T$, *i.e.*, the domain $\Omega^{\text{ref}}$ in Figure 2 is the image of the mapping $T$:

$$\Omega^{\text{ref}} = T[\tilde{\Omega}], \tag{50}$$

where $\tilde{\Omega}$ is the parameter space of $T$. For analysis purposes, we only consider analysis-suitable spline model made of multiple patches with matching interfaces. The spline basis functions of $T$ are used to defined the finite element spaces (the $N_i$ in Equation (6)). Additionally, the mappings $G^{(s)}$ are given as Bezier multivariates, which might be obtained via a Bezier decomposition of a spline mesh. Finally, we combine these two mappings to build an explicit representation of the lattice cells. This is done via functional composition, meaning that the subdomains are the images of the composed mappings:

$$\Omega^{(s)} = (G^{(s)} \circ T)[\tilde{\Omega}]. \tag{51}$$

More information regarding this multiscale geometric representation can be found in [16, 28] for instance.

**Remark 8.** *Despite the fact that we leverage the use of spline functions for both, the geometry description and the solution discretization, following the isogeometric paradigm, the method presented is completely agnostic to type of finite element discretization applied to the PDE.*

### 4.2. A multiscale Greedy procedure

#### 4.2.1. Assembly via lookup tables

As quickly mentioned previously, in order to build the reduced basis of the local stiffness matrices as introduced in Equation (38), we employ a greedy approach. This enables to get a reduced basis that interpolates some local stiffness matrices, *i.e.*, those forming the principal sub-problems. The crucial point here is that we do not build the reduced basis on the local stiffness matrices directly. This might be too expensive. Instead, the similarities between the cells can be observed by looking at the weak formulations given in Equation (5). After a change of coordinates by pulling back the macro-mappings $G^{(s)}$, all local bilinear forms are defined over the reference domain $\Omega^{\text{ref}}$, as:

$$a^{(s)}(u, v) = \int_{\Omega^{\text{ref}}} \hat{\varepsilon}(u) : \hat{\mathbf{C}}^{(s)} : \hat{\varepsilon}(v) \, d\xi, \tag{52}$$

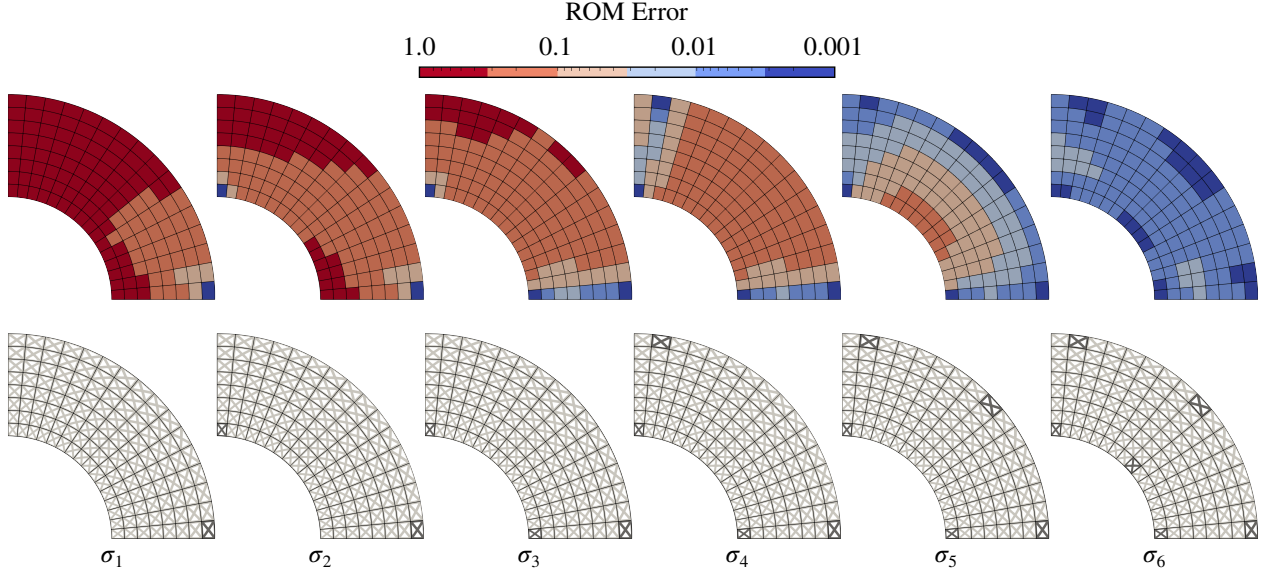

Figure 5: Iterative identification of the principal subdomains via the greedy procedure occurring at the macro-scale. The ROM Error corresponds to the norm of the residuals $\|\mathbf{\Delta}_i^{(s)}\|_\infty$ as written in Algorithm 4.

where $\hat{\mathbf{C}}^{(s)} : \Omega^{\text{ref}} \to \mathbb{R}^{3\times3\times3\times3}$ are fourth-order tensor fields which encapsulate all the dependencies on the macro-mappings. I.e., it encodes the material properties of the tensors $\mathbf{C}^{(s)}$ as well as all the quantities depending on the macro-mappings $G^{(s)}$. So, one can take a look to these fields in order to build the reduced basis, instead of doing it over the local stiffness matrices. Furthermore, as all the other operators and quantities in Equation (52) are the same for every cell, one may resort to the approach developed in [28]: It is possible to build lookup tables with precomputed integrals in order to assemble quickly local stiffness matrices. This approach consists in introducing suitable polynomial approximations of the tensor fields:

$$\hat{\mathbf{C}}^{H(s)}(\xi) = \mathbf{N}^\top(\xi)\mathbf{A}^{(s)}, \quad i.e., \quad \hat{C}^{H(s)}_{ijkl}(\xi) = \sum_{q=1}^{n_A} N_q(\xi) A^{(s)}_{ijkl,q}, \tag{53}$$

where $\mathbf{N}$ contains $n_A$ polynomial basis functions and $\mathbf{A}^{(s)}$ collects the associated polynomial coefficients. The latter ones are constants and encapsulate the dependencies on the macro-mappings $G^{(s)}$ and the material properties. Consequently, one can take them out of the integrals when substituting Equation (53) into Equation (52), which leads to integrals that are independent of the macro-mappings and shared by all the cells. Those integrals are defined over the reference domain $\Omega^{\text{ref}}$, computed once and for all, and stored in lookup tables $\mathbf{T}^{\text{ref}}$. Building the stiffness matrices is then done by combining those precomputed integrals together with the polynomial coefficients, *i.e.*, the assembly reads as:

$$\mathbf{K}^{(s)} = \mathbf{T}^{\text{ref}} \otimes \mathbf{A}^{(s)}, \tag{54}$$

where the operator $\otimes$ encodes multiplications, summations, and reordering of both quantities. We refer again to [28] where this multiscale assembly strategy is detailed.

**Remark 9.** *It should be noticed that applying the fast assembly strategy via lookup tables [28] is just a choice. The presented solver does not depend on this method. It however offers additional advantages as being able to perform matrix-vector products between the $\mathbf{K}^{(s)}$ and input vectors in a matrix-free fashion (i.e., without explicitly building the matrices, but by using the lookup table and the polynomial coefficients directly, to save memory for instance). As will be presented in Section 4.2.2, it also simplifies the construction of the reduced basis of the stiffness matrices as one can operate on a small collections of coefficients $\mathbf{A}^{(s)}$ instead of functions $\hat{\mathbf{C}}^{(s)}(\xi)$ or large data $\mathbf{K}^{(s)}$.*

**Algorithm 4** Greedy algorithm for the identification of the principal subdomains.

| | |
|---|---|
| Initialization | |
| $\forall s,\ \text{compute}:\ \mathbf{\Delta}_0^{(s)} = \mathbf{A}^{(s)}/\|\mathbf{A}^{(s)}\|_2$ | |
| For $i = 1, 2, \ldots$, until convergence | |
| $\sigma_i = \underset{s=1,\ldots,N}{\mathrm{argmax}} \|\mathbf{\Delta}_{i-1}^{(s)}\|_\infty$ | Find current principal subdomain |
| $\boldsymbol{\zeta}_i = \mathbf{\Delta}_{i-1}^{(\sigma_i)}/\|\mathbf{\Delta}_{i-1}^{(\sigma_i)}\|_2$ | Add new vector basis |
| $\forall s,\ \text{compute}:\ \beta_i^{(s)} = \mathbf{\Delta}_{i-1}^{(s)} \cdot \boldsymbol{\zeta}_i$ | Compute additional coordinates |
| $\forall s,\ \text{compute}:\ \mathbf{\Delta}_i^{(s)} = \mathbf{\Delta}_{i-1}^{(s)} - \beta_i^{(s)} \boldsymbol{\zeta}_i$ | Update residuals |
| $\underset{s=1,\ldots,N}{\max} \|\mathbf{\Delta}_i^{(s)}\|_\infty \overset{?}{\lessgtr} tol_{\text{RB}}$ | Check convergence |
| Outputs: | |
| $\mathbb{Z} = [\boldsymbol{\zeta}_1 \ldots \boldsymbol{\zeta}_{N_{\text{rb}}}]$ | Orthonormal basis |
| $\mathrm{S} = \{\sigma_1, \ldots, \sigma_{N_{\text{rb}}}\}$ | Principal subdomain indices |
| $\mathring{\mathbf{A}}^{(s)} = \|\mathbf{A}^{(s)}\|_2 \mathbb{Z} \boldsymbol{\beta}^{(s)},\ \forall s$ | Certified affine decompositions |

#### 4.2.2. Algorithm to identify the principal subdomains

In this work, the procedure to identify the principal subdomains relies on a greedy algorithm (which can be viewed as a truncated Gram-Schmidt process), which is a common technique to build reduced bases, see for instance [25, 54]. As discussed in Remark 9, the snapshot matrix could be formed using either the stiffness matrices $\mathbf{K}^{(s)}$, the functions $\hat{\mathbb{C}}^{(s)}(\xi)$ in the bilinear forms (see Equation (52)), or the polynomial coefficients $\mathbf{A}^{(s)}$ as defined by Equation (53). Using directly the stiffness matrices is obviously the more expensive as it contains a lot of data. It is more efficient to consider one of the two other options. Here, we work directly with the polynomial coefficients (that come from the polynomial approximations of the fourth-order tensor fields $\hat{\mathbb{C}}^{(s)}$). The greedy process for selecting the principal subdomains and performing the affine decomposition is given in Algorithm 4, and illustrated in Figure 5. Interestingly, due to its construction, the reduced space interpolates the principal subdomains, *i.e.*:

$$\mathring{\mathbf{A}}^{(\sigma_i)} = \mathbf{A}^{(\sigma_i)}, \quad \forall i = 1, \ldots, N_{\text{rb}}. \tag{55}$$

Consequently, we can perform a change of basis from the orthonormal basis $\mathbb{Z}$ formed during the greedy process to the basis formed by quantities defined over the principal subdomains $\mathbb{A} = [\mathbf{A}^{(\sigma_1)}, \ldots, \mathbf{A}^{(\sigma_{N_{\text{rb}}})}]$, *i.e.*:

$$\mathring{\mathbf{A}}^{(s)} = \|\mathbf{A}^{(s)}\|_2 \mathbb{Z} \boldsymbol{\beta}^{(s)} = \mathbb{A} \boldsymbol{\alpha}^{(s)}. \tag{56}$$

As a direct consequence, we end up with the approximate affine decompositions of the stiffness matrices as given in Equation (38) which is the starting point of our inexact FETI-DP solver.

Instead of using the polynomial coefficients $\mathbf{A}^{(s)}$, one could perform the affine expansion of the functions $\hat{\mathbb{C}}^{(s)}$ via the Empirical Interpolation Method. It consists in a similar greedy process as the one defined in Algorithm 4, but it involves an additional step where so-called magic points are selected. The interested readers may refer to [25, 54] for more details regarding EIM and reduced basis space construction in general.

## 5. Numerical examples

All the following numerical examples were computed on a 2.50GHz i7-11850H (11th Gen Intel(R) Core(TM)) processor with 32 GB RAM (serial computing only). More precisely, three 2D examples and two 3D examples are investigated in this section. They all share same material parameters, *i.e.*, an isotropic elastic material with Young's modulus $E = 5000$ MPa and Poisson's ratio $\nu = 0.40$. Regarding the solver, the global convergence relative tolerance was set to $tol_{\text{GMRES}} = 1e-5$, while the local convergence relative tolerance for the interface problem was set to $tol_{\text{CG}} = 1e-11$. The implementation has been done with Python using numpy/scipy packages [23, 67]. It also calls Cholesky factorization routines provided by scikit-sparse package which makes use of CHOLMOD from the SuiteSparse linear algebra package [11].

### 5.1. 2D examples

#### 5.1.1. Description

Let us start with the two-dimensional examples. The first inputs were the macro-geometries which are depicted in Figure 6. All these geometries were geometrically described via spline based models. They are of three kinds:

1. Fig. 6(a) – A rectangular domain subjected to an imposed displacement. The geometry was described with a linear B-Spline patch.

2. Fig. 6(b) – A curved beam under bending. The geometry was (exactly) described with a quadratic NURBS patch.

3. Fig. 6(c) – A brake pedal where the loading case mimics the operating conditions. The geometry was also described with a quadratic NURBS patch with knot-vectors $U = \{0, 0, 0, 0.5, 0.5, 1, 1, 1\}$ and $V = \{0, 0, 0, 1, 1, 1\}$. The coordinates and weights of the control points are given in the mentioned figure.

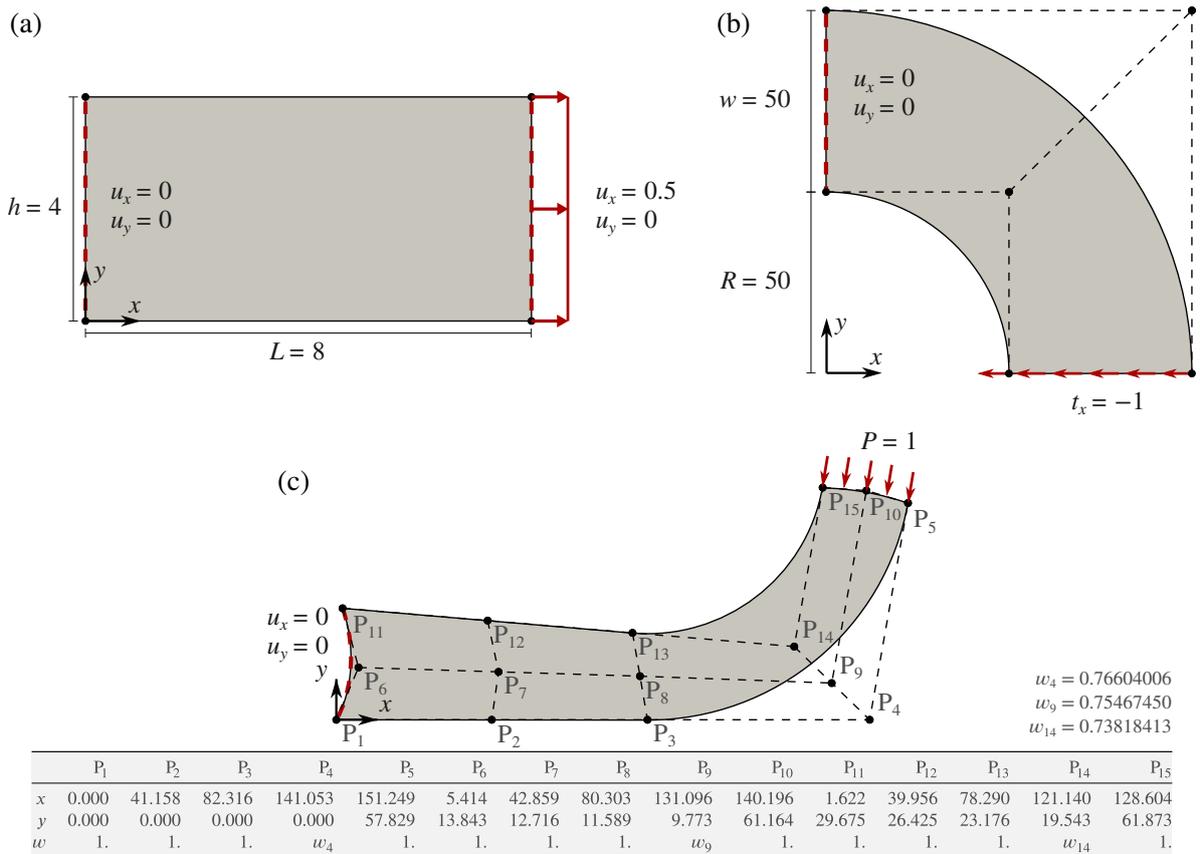

|   | $P_1$ | $P_2$ | $P_3$ | $P_4$ | $P_5$ | $P_6$ | $P_7$ | $P_8$ | $P_9$ | $P_{10}$ | $P_{11}$ | $P_{12}$ | $P_{13}$ | $P_{14}$ | $P_{15}$ |
|---|---|---|---|---|---|---|---|---|---|---|---|---|---|---|---|
| $x$ | 0.000 | 41.158 | 82.316 | 141.053 | 151.249 | 5.414 | 42.859 | 80.303 | 131.096 | 140.196 | 1.622 | 39.956 | 78.290 | 121.140 | 128.604 |
| $y$ | 0.000 | 0.000 | 0.000 | 0.000 | 57.829 | 13.843 | 12.716 | 11.589 | 9.773 | 61.164 | 29.675 | 26.425 | 23.176 | 19.543 | 61.873 |
| $w$ | 1. | 1. | 1. | $w_4$ | 1. | 1. | 1. | 1. | $w_9$ | 1. | 1. | 1. | 1. | $w_{14}$ | 1. |

Figure 6: The macro-geometries and the loading scenarios that define the 2D examples: (a) a rectangular domain subjected to an imposed displacement, (b) a curved beam under bending, and (c) a brake pedal subjected to operating conditions.

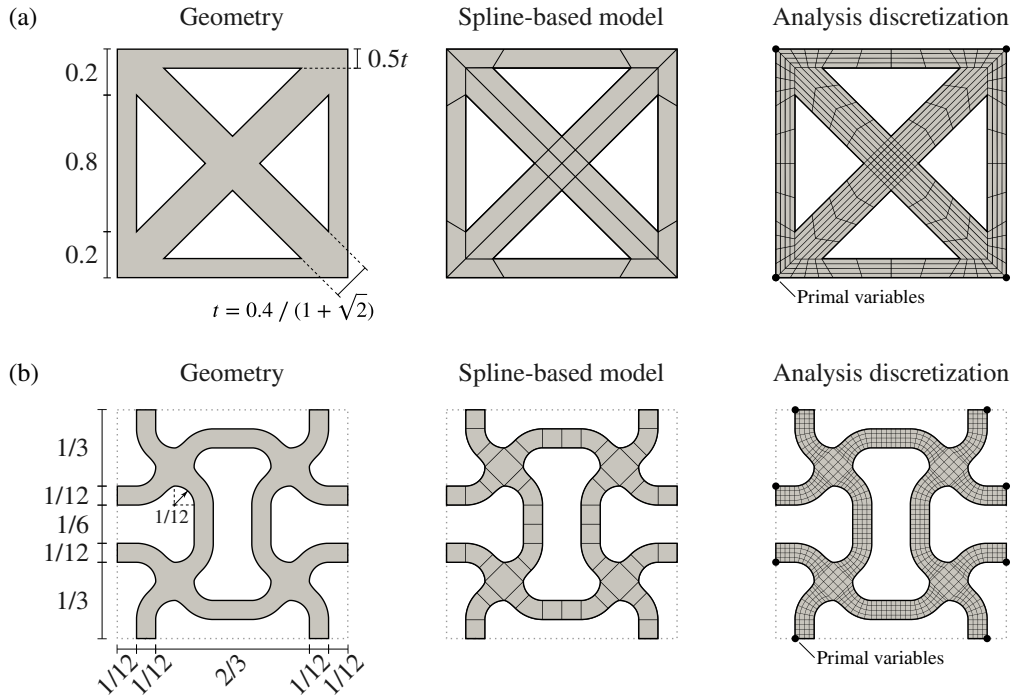

Figure 7: The 2D cells used during the numerical studies: (a) a hollow square with a cross along its diagonal, and (b) an auxetic lattice pattern.

Then, we associated a reference cell (microscale) to each of these macro-geometries. Again, we modelled these cells with spline based models. Two cells were considered as depicted in Figure 7:

1. Fig. 7(a) – A hollow square linked by a cross. It required twenty-four linear Bézier patches to build an analysis-suitable model (*i.e.*, matching interfaces between the patches, and no trimmed patches, see, *e.g.*, [8]). The model was further refined to generate the analysis model. This analysis discretization was made of 384 cubic elements, which led to 1866 DOF per cell.

2. Fig. 7(b) – An auxetic lattice pattern (obtained through a topology optimization procedure in [68]). We modelled the geometry with twenty quadratic matching NURBS patches. The refined version used during the analysis was made of 832 cubic elements for a total of 4368 DOF per cell.

We combined the crossed hollow square unit cell (Figure 7(a)) with the curved beam problem (Figure 6(b)) and the brake pedal problem (Figure 6(c)). The auxetic lattice (Figure 7(b)) was combined with the rectangular macro-domain subjected to an imposed displacement (Figure 6(a)).

Several results from these structural analyses are shown in Figure 8. We depict the displacement field and the von Mises stress field for the three test cases. As we perform full-scale simulations, we directly have access to these quantities over the different cell-struts.

### 5.1.2. Solver performance

We study the performance of the developed solver on these 2D numerical examples. A first question concerns the influence of the reduced bases. For instance, how many principal cells is it suitable to pick? What is the influence of the size of the reduced bases on the iterative convergence of the solver?

*Size of the reduced bases.* Figure 9 gives an insight of the behavior of the solver *w.r.t.* the number of principal cells, (*i.e.*, by varying on the tolerance $tol_{\text{RB}}$ during the greedy process, see Table 4). We performed this study on the curved beam problem (Figure 8(b)). More precisely, the corresponding results read as follows:

- 27 iterations when 6 principal cells out of the 512 cells were selected ($tol_{\text{RB}} = 1e - 1$),

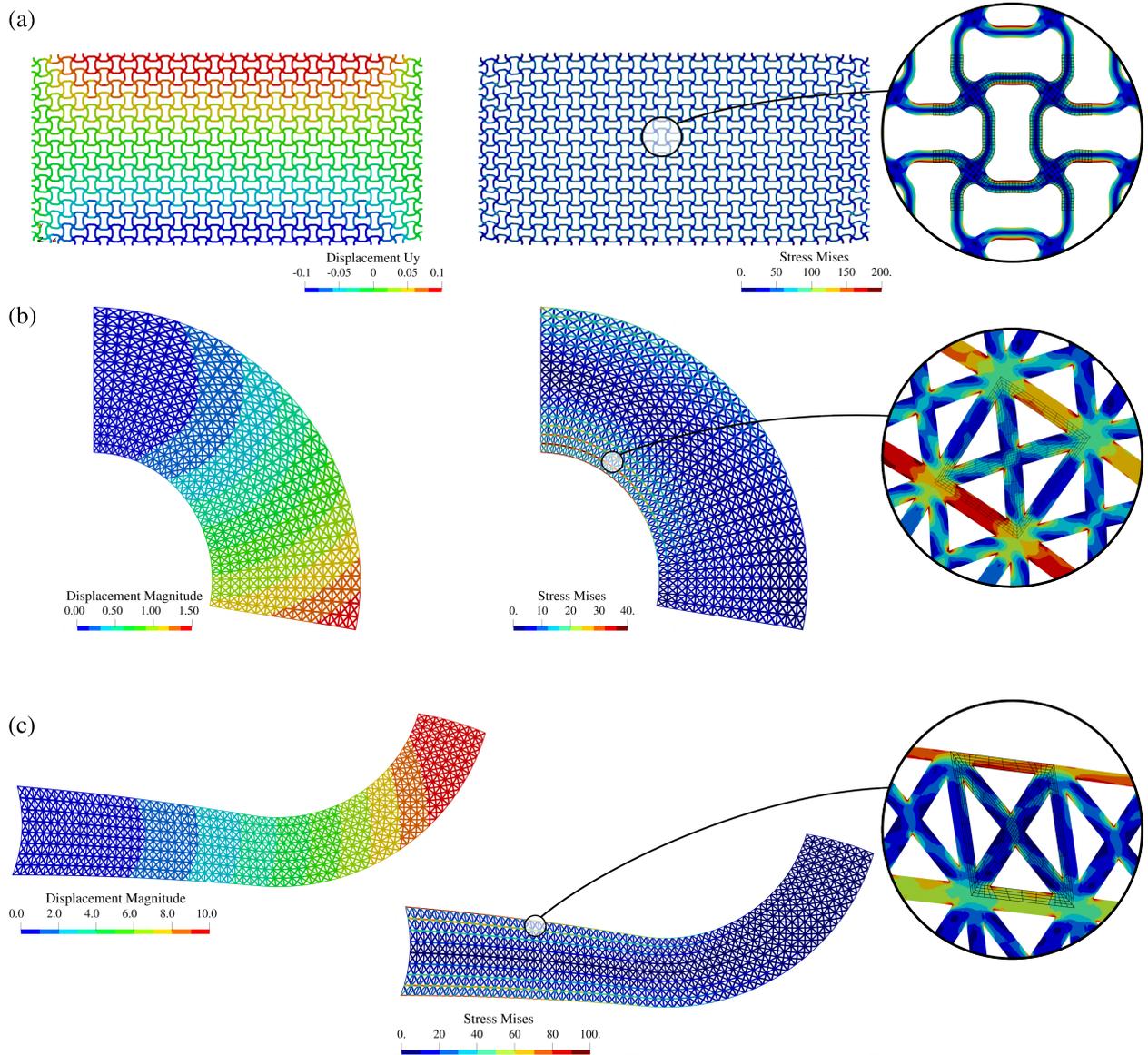

Figure 8: Results of the structural analyses for the 2D test cases. Note that the auxetic behavior can be observed in (a): The initial macro rectangle inflates while it was subjected to traction.

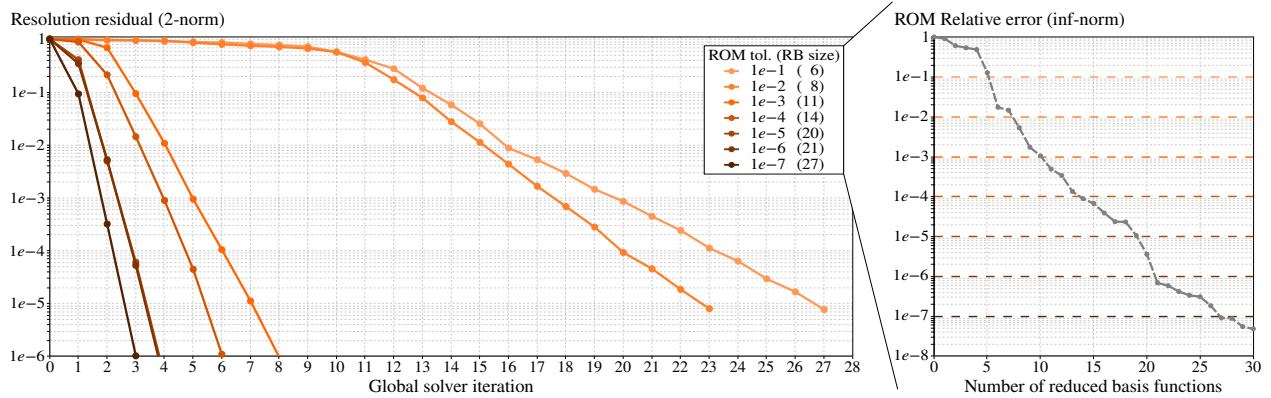

Figure 9: Influence of the reduced basis on the convergence of the solver. The larger the reduced basis, the lower the number of iterations the solver needs.

- 23 iterations when 8 principal cells out of the 512 cells were selected ($tol_{RB} = 1e-2$),
- 8 iterations when 11 principal cells out of the 512 cells were selected ($tol_{RB} = 1e-3$),
- 6 iterations when 14 principal cells out of the 512 cells were selected ($tol_{RB} = 1e-4$),
- 4 iterations when 20 principal cells out of the 512 cells were selected ($tol_{RB} = 1e-5$),
- 4 iterations when 21 principal cells out of the 512 cells were selected ($tol_{RB} = 1e-6$),
- 3 iterations when 27 principal cells out of the 512 cells were selected ($tol_{RB} = 1e-7$).

This observed behavior was to be expected: the larger the RB (*i.e.*, the lower $tol_{RB}$), the faster the convergence (in terms of number of global iterations). In terms of computational cost, it is straightforward however that increasing the number of principal cells leads to more memory usage. Indeed, it leads to compute and to store more quantities, as for instance, the Cholesky factorizations of the principal local problems. Regarding the computational time, selecting a smaller number of principal cells decreases the time needed to build the ROMs, and the time per global iteration. However, as it leads to a larger number of global iterations, the overall computational time might be above the one obtained with a larger number of principal cells. And this is somehow what we observed in the curved beam example. Indeed, the overall computational time was:

- around 25 s for the cases with 6 and 8 principal cells (respectively, $tol_{RB} = 1e-1$ and $tol_{RB} = 1e-2$),
- around 11 s for all the other cases, with a small decrease while increasing the RB size (*i.e.*, lower tolerance while building the RBs in the greedy process).

We experienced the same behavior for other test cases: A larger number of principal cells leads to fewer global iterations, it however increases (a bit) the memory usage, but it leads to equivalent or less overall computational time. Thus, we advise to set a relative tolerance in the greedy selection algorithm (see again Table 4) between $1e-5$ and $1e-7$. This tolerance might seem low, but one should keep in mind that this tolerance is set over the local stiffness matrices only, and not over the local DD operators directly (*e.g.*, the primal and dual Schur complements, see again Equation (39) and Remark 2). The particularity here is that we cannot build the ROMs *a priori*, and even more, we do not have access to the quantities for which we are building ROM-based approximations. For instance, we cannot run a similar greedy process on the Schur complements, directly, as it would mean to build them all, which would take more time than simply solving the problem with a standard FETI-DP approach. That is why, by selecting a sufficiently large number of principal cells, we ensure that the ROMs of the local DD operators are suitable.

**Remark 10.** *The macro-geometry of the curved beam is modelled with a NURBS, see Figure 6(b). Hence, the quarter circles are exactly represented but not parameterized in the sense of the arc length. Thus, the cells do not follow exactly a circular repeating pattern. Furthermore, we do not apply a change of basis to rotate the cells prior the identification of the principal cells. Using an algebraic macro-mapping together with a change of basis could reduce the number of principal cells for this particular example.*

| Test case | | | ROM-based inexact FETI-DP | | | | Standard FETI-DP | | |
|---|---|---|---|---|---|---|---|---|---|
| *Macro-geom.* | *#Cells* | *#DOFs* | *Memory* | *Comp. Time* | *Glo. It.* | *Loc. It.* | *Memory* | *Comp. Time* | *Loc. It.* |
| Rectangle | 16 × 08 | 0.5M | 0.0 GiB | 0.3s (0s, 0s) | 1 | 28 | 1.6 GiB | 8s ( 4s, 4s) | 28 |
| Rectangle | 32 × 16 | 2.2M | 0.1 GiB | 1s (0s, 1s) | 1 | 29 | 6.4 GiB | 38s (25s, 14s) | 29 |
| Rectangle | 64 × 32 | 8.9M | 0.2 GiB | 5s (1s, 4s) | 1 | 29 | > available | - | - |
| Curved beam | 16 × 08 | 0.2M | 0.2 GiB | 5s (3s, 2s) | 3 | 35 | 0.8 GiB | 6s ( 4s, 2s) | 35 |
| Curved beam | 32 × 16 | 1.0M | 0.3 GiB | 11s (3s, 8s) | 3 | 35 | 3.2 GiB | 22s (13s, 9s) | 35 |
| Curved beam | 64 × 32 | 3.8M | 0.5 GiB | 38s (3s, 35s) | 3 | 35 | 13.4 GiB | 73s (38s, 35s) | 34 |
| Brake pedal | 32 × 04 | 0.2M | 0.3 GiB | 10s (6s, 4s) | 3 | 42 | 0.8 GiB | 6s ( 3s, 3s) | 42 |
| Brake pedal | 64 × 08 | 1.0M | 0.3 GiB | 21s (6s, 15s) | 4 | 43 | 3.2 GiB | 25s (14s, 11s) | 42 |
| Brake pedal | 128 × 16 | 3.8M | 0.5 GiB | 86s (5s, 81s) | 5 | 43 | 12.7 GiB | 80s (38s, 42s) | 43 |

Table 1: Performance of the solver in the case of 2D problems.

*Scalability.* In order to further investigate the performance of our ROM-based inexact FETI solver, we performed a scalability test for the three 2D examples. The results are presented in Table 1. The table also provides the results obtained with a standard FETI-DP solver, which means that the interface problem is "exactly" formed. In practice, it means that all the local stiffness matrices have been factorized prior to running the iterative algorithm. This is why, in the *Memory* column for the Standard FETI-DP data in Table 1, one can notice a substantial amount of memory usage in comparison with the *Memory* usage involved when employing the developed approach where we limit the construction of local operators to the picked principal cells. This is a great advantage of the proposed approach. Despite the use of full-scale models of lattice structures, the memory requirements were limited as we exploit the inherent repetitive pattern of those structures. On the contrary, black-box solvers are not aware of the structure of the problem at hand, and, hence, the available computational resources were rapidly surpassed by these (standard) approaches. This point is supported by the results: We were not able to compute the most refined version of the auxetic rectangular problem (64-by-32 cells for a total of 8.9M) with the standard FETI-DP approach and our computational resources, while we easily succeeded with the developed ROM-based FETI-DP approach.

The number of iterations required to solve the interface problems was the same with both approaches (*Loc. It.* columns). This was expected as the FETI-DP is known to be scalable for 2D problems when the so-called corners are chosen as the primal variables, see for instance [37]. The increase of the number of total cells has little influence on the number of global iterations for the ROM-based solver, as observed, for instance, for the brake pedal problem. However, this is minor as the solver still converged in very few global iterations (5 in the finest configuration). Furthermore, Table 1 provides in the columns *Comp. Time* the computational time in the following format:

- *total time* (*time during preprocessing step*, *time for the iterative algorithm*).

For the proposed ROM-based FETI-DP, the preprocessing step consists in building the reduced models $\mathring{(\cdot)}^{(s)}$ that approximate the DD local operators (see, *e.g.*, Equation (43)), whereas in the standard FETI-DP, it consists mainly in factorizing all the local operators. Our test cases have an increasing level of complexity, in the sense of the geometry. For the rectangular problem, all the cells were actually the same. For the curved beam, they exhibited rather small geometric differences, whereas for the brake pedal, each cell had a different shape. Regarding the first example, the presented approach shows significant speedup in comparison with the standard approach. Even in the finest configuration with 8.9M DOF, the solver took only 5 s to get the solution.

**Remark 11.** *The outstanding solver performance for the rectangular domain example is due to the fact that all the cells are geometrically identical. In this particular case, the ROM introduces no approximations, being the reduced bases composed of a single basis entry and all DD operators exactly recovered. Consequently, only one single global iteration is required. Indeed, the role of the global iterations in our solver is to compensate the ROM induced errors. If no approximations are introduced by the ROM, as it is the case for this particular example, then the inexact FETI-DP preconditioner is nothing else than the inverse of the system matrix, i.e., $\mathring{\mathbf{A}}^{-1} = \mathbf{A}^{-1}$.*

For the curved beam, our solver took half the time needed with the standard FETI-DP, and for the brake pedal, it took a similar amount of time. For this last case, the number of principal cells was about 30 which explain the increase in computational time. However, as said previously, the memory usage was much less with the ROM-based approach. Applying blindly standard FETI-DP approach for large lattice structure would be unreasonable for this reason. This will become even clearer in the following 3D examples. We also think that some speedup with the ROM-based solver could be further obtained. Indeed, the code is partly built upon the sparse library from scipy which does not provide top performing matrix-matrix and matrix-vector products (contrary to the dense linear algebra which makes use of BLAS). Approximating the local solutions via the RB-based process described in Equation (46) took the major part of the time.

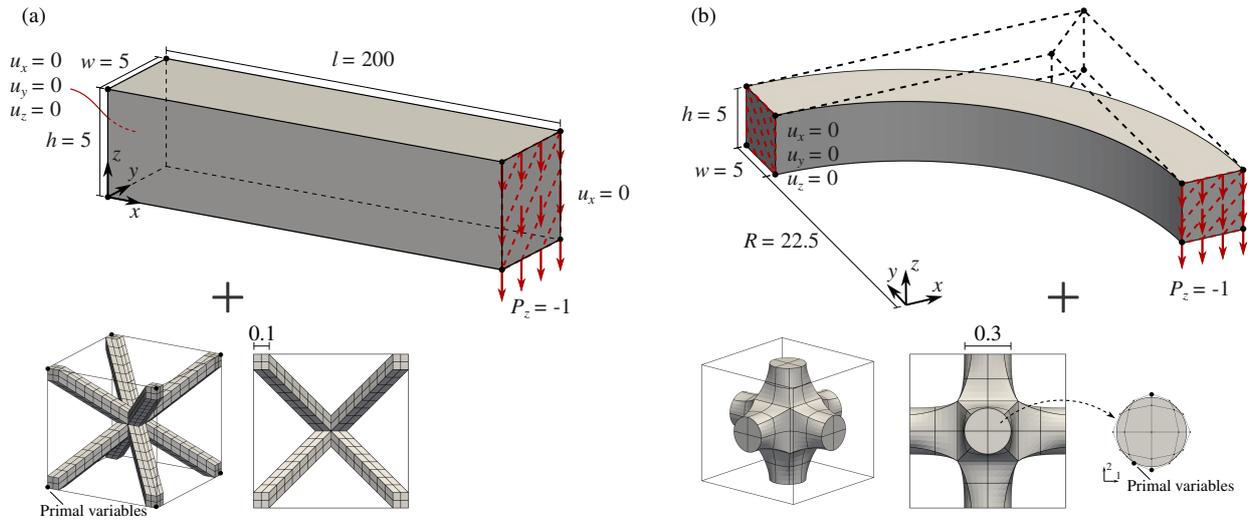

Figure 10: The macro-geometries, the loading scenarios, and the reference cells that define the 3D examples: (a) a straight beam filled by Body-Centered Cubic patterns, and (b) a curved beam filled by patterns initiated from a reference cell made of a cross with circular arms.

### 5.2. 3D examples

#### 5.2.1. Description

We now study two 3D examples which are depicted in Figure 10. The first test case consists in a straight beam under bending. Due to the symmetry of the problem, we considered only half of the problem (we could even consider only one quarter of the problem). The macro-geometry was modeled with a linear B-Spline patch with one single element. We associated to this macro-geometry a reference lattice cell. In this case, it was a Body-Centered Cubic cell. This BCC cell was modeled with 32 linear B-Spline patches. The geometric model was further refined for analysis purposes. The analysis discretization associated to one cell included 832 quadratic spline elements for a total of 9021 DOF. Several results from a structural analysis are given in Figure 11. Again, as we solved the full-scale problem, we had direct access to the displacement and the stress field over the different cell-struts.

The second test case (see again Figure 10(b)) consists in a curved beam which was subjected to bending too. The macro-geometry was modeled with a quadratic NURBS with one element. This enabled to exactly represent the circular shape of the beam. We associated to this macro-geometry a reference lattice cell whose shape looks like a cross with curved and circular arms. Initially, the geometric model of the reference cell was made of 7 cubic patches which were then further refined such that the analysis model included 56 cubic elements, and a total of 2175 DOF. We provide an illustration of the displacement and stress field for this test case in Figure 12.

#### 5.2.2. Solver performances

We performed similar experiments as for the 2D examples. More specifically, we study the behavior of the presented ROM-based FETI-DP solver while increasing the number of cells for each problem. The results are summarized in Table 2, that presents the same structure as Table 1.

As already mentioned, an important characteristic of the proposed solver is the drastic reduction of the memory usage in comparison to "black box" solvers, as for instance a standard FETI-DP solver. This has been seen on the 2D examples, and it is even more obvious for the 3D examples. Here, only the coarser configurations could be solved with the standard method. Indeed, the memory limit was rapidly reached with the classic FETI-DP, whereas the developed ROM-based FETI-DP solver enabled us to simulate geometries with a higher number of lattice cells and low memory resources. Regarding the results presented in Table 2, enormous speedup was obtained for the straight beam problem. The reasons are exactly the same as those mentioned for the first of the 2D examples, *i.e.*, the rectangle with the auxetic lattice pattern. Here, there were, by construction, no approximations of the local operators, and the reduced bases included one single element only as all the cells were the same. Thus, the solver took less than one minute to solve a problem with $6 \times 6 \times 24$ cells and 7.8M DOF. Regarding the curved beam problem, it required 3 or 4 global iterations to converge. It is faster than the standard FETI-DP, and required much less memory. Thus, it was possible to run on a simple laptop a case with $6 \times 6 \times 56$ cells and 3.7M DOF with the ROM-based FETI-DP solver. Indeed, only 20 principal cells were selected out of the 2016 cells. On the contrary, the standard approach reached the memory limit much before, and only a case with $4 \times 4 \times 32$ cells and 1.1M DOF was trackable.

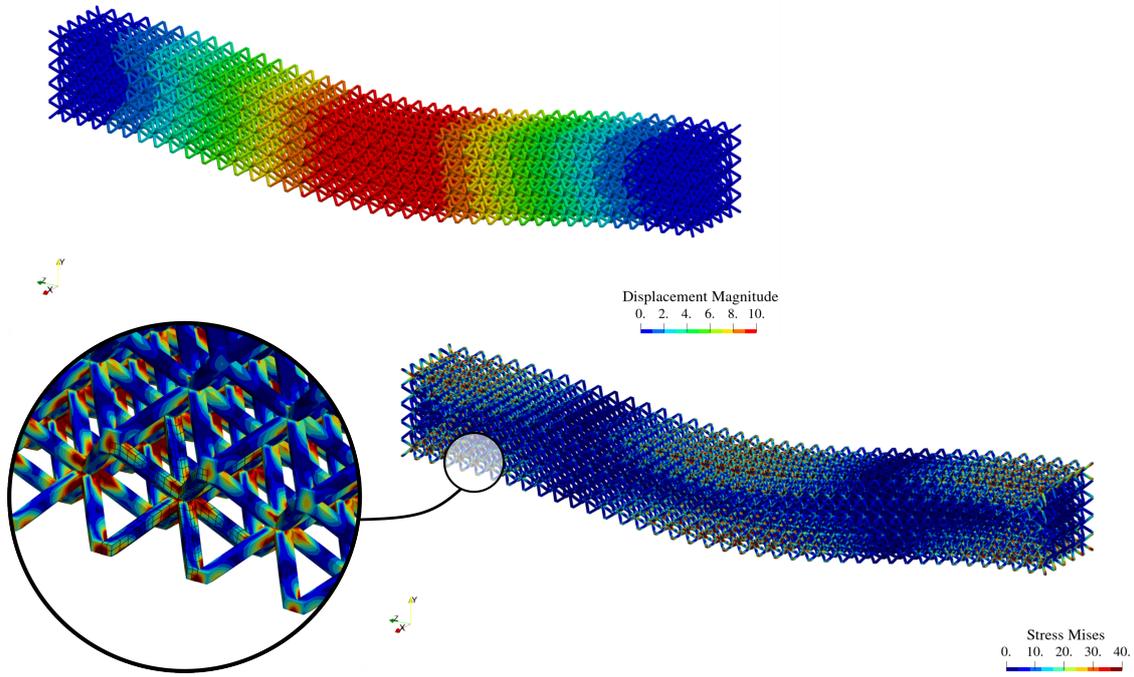

Figure 11: Structural analysis results for the straight 3D beam problem.

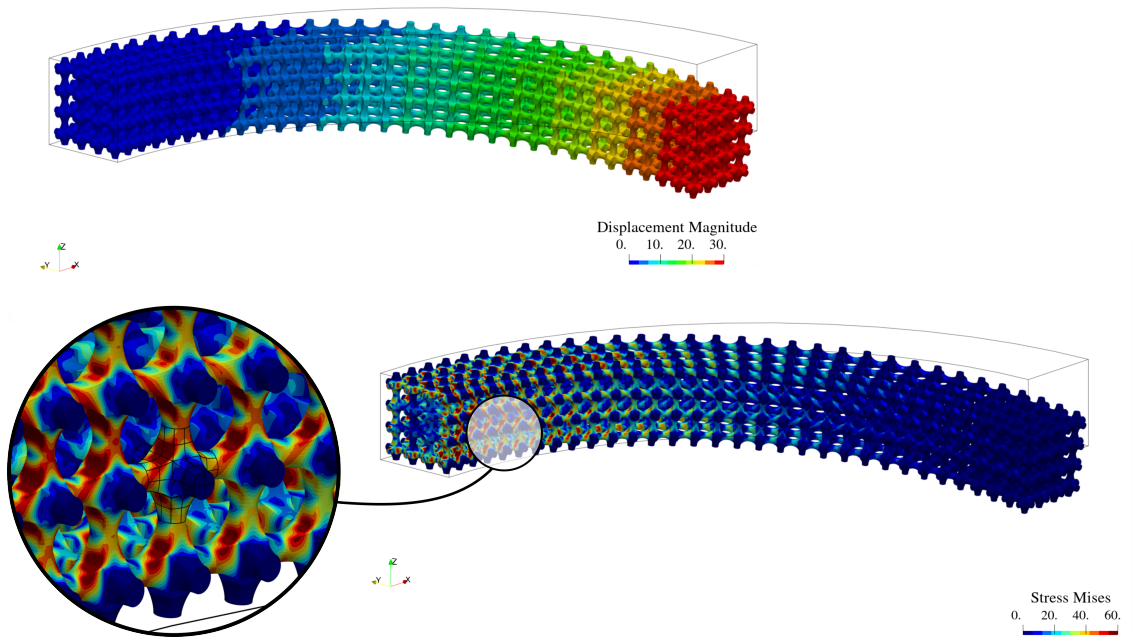

Figure 12: Structural analysis results for the curved 3D beam problem.

| Test case | | | ROM-based inexact FETI-DP | | | | Standard FETI-DP | | |
|---|---|---|---|---|---|---|---|---|---|
| *Macro-geom.* | *#Cells* | *#DOFs* | Memory | Comp. Time | Glo. It. | Loc. It. | Memory | Comp. Time | Loc. It. |
| Straight beam | $2 \times 2 \times 08$ | 0.3M | 0.5 GiB | 0min 3s (1s, 2s) | 1 | 103 | 4.2 GiB | 0min 36s (10s, 26s) | 103 |
| Straight beam | $3 \times 3 \times 12$ | 1.0M | 0.9 GiB | 0min 7s (2s, 5s) | 1 | 106 | 14.0 GiB | 2min 4s (31s, 93s) | 106 |
| Straight beam | $4 \times 4 \times 16$ | 2.3M | 1.7 GiB | 0min 14s (3s, 11s) | 1 | 108 | > available | - | - |
| Straight beam | $5 \times 5 \times 20$ | 4.5M | 3.1 GiB | 0min 26s (4s, 22s) | 1 | 109 | > available | - | - |
| Straight beam | $6 \times 6 \times 24$ | 7.8M | 5.1 GiB | 0min 44s (6s, 38s) | 1 | 110 | > available | - | - |
| Curved beam | $3 \times 3 \times 24$ | 0.5M | 1.2 GiB | 0min 36s (8s, 24s) | 4 | 93 | 6.9 GiB | 0min 49s (23s, 26s) | 92 |
| Curved beam | $4 \times 4 \times 32$ | 1.1M | 1.6 GiB | 1min 0s (12s, 48s) | 3 | 103 | 15.7 GiB | 1min 53s (49s, 64s) | 102 |
| Curved beam | $5 \times 5 \times 40$ | 2.2M | 2.2 GiB | 1min 57s (19s, 98s) | 3 | 113 | > available | - | - |
| Curved beam | $6 \times 6 \times 56$ | 3.7M | 3.1 GiB | 3min 40s (20s, 230s) | 4 | 120 | > available | - | - |

Table 2: Performance of the solver in the case of 3D problems.

## 6. Conclusions

In this work, we developed a novel method for performing full fine-scale analysis of architected materials at very low computational cost. This dedicated solver takes advantage of the intrinsic geometrical nature of such structures; that is, the non-overlapping domain decomposition into cells that are all similar. The faced challenge was: Is it possible to leverage on those similarities to save memory usage and computational time during the simulation? Resorting to the class of DD methods, especially to the inexact variants, and complementing the approach with ROM, this turned out possible. More precisely, we built upon inexact FETI-DP algorithms that offer the possibility of not solving exactly all the local, cell-wise systems. Then, to design an efficient block preconditioner, we proposed to extract the principal cell-wise stiffness matrices with a greedy approach, and to use them to approximate all the local quantities occurring in a standard FETI-DP algorithm. This results in a solver which only needs, in terms of local solutions, the ones involving the principal stiffness matrices. Therefore, numerous local factorizations are saved compared to standard DD algorithms. Furthermore, although not mandatory in our methodology, we adopted a geometrical modeling based on spline composition and performed, following the isogeometric paradigm, analysis of lattices using look-up tables and macro-fields (see [28]). This allowed us to extract the principal local operators by acting only on the macro-fields, which avoids forming all the local operators, and therefore makes our method also multi-scale and matrix-free.

Several numerical experiments in 2D and 3D and covering different micro- and macro-mappings were performed to assess the solver performance. Only serial computations were conducted at this stage. Comparisons with the standard FETI-DP were also carried out to highlight the interest of our solver. In any case, the memory usage was drastically reduced with our solver compared to the standard FETI-DP, which allowed us to compute problems of several millions of DOF within few minutes using an off-the-shelf laptop. Eventually, the computational time was also very low, in particular when the macro-mapping is affine (since all the cells become identical), which is usually the case in real architected material designs. The obtained results suggest that introducing multiscale modeling for lattice structures to keep computing resources reasonable could be circumvented with such HPC algorithms. In some sense, the multiscale, heterogeneous, and periodicity features of the structure are used directly in the solver, so we do not need to bring them at the modeling stage. In this respect, an efficient parallel implementation, and the extension of the method to the nonlinear regime, as well as to take into account manufacturing-induced defects, surely constitute interesting prospects to this work.

*Acknowledgment.* This research was supported by the European Union Horizon 2020 research and innovation program, under grant agreement No 862025 (ADAM2) and by the Swiss National Science Foundation through the project No 40B2-0_187094 (BRIDGE Discovery 2019). The second author, Robin Bouclier, also gratefully acknowledges support from the French National Research Agency under Grant ANR-22-CE46-0007 (AVATAR).